\documentclass[a4paper,12pt]{amsart}

\usepackage{amssymb}
\usepackage{yfonts}

\def\hbar{\bar{h}}

\def\iso{\buildrel \sim\over\to}

\def\GA{{\mathfrak{A}}}

\def\GS{{\mathfrak{S}}}

\def\Gh{{\mathfrak{h}}}

\def\Gsl{{\mathfrak{sl}}}

\def\CA{{\mathcal{A}}}
\def\CB{{\mathcal{B}}}
\def\CC{{\mathcal{C}}}
\def\CD{{\mathcal{D}}}

\def\CF{{\mathcal{F}}}

\def\CH{{\mathcal{H}}}

\def\CL{{\mathcal{L}}}

\def\CO{{\mathcal{O}}}
\def\CP{{\mathcal{P}}}

\def\CR{{\mathcal{R}}}

\def\CU{{\mathcal{U}}}

\def\BA{{\mathbf{A}}}
\def\BB{{\mathbf{B}}}
\def\BC{{\mathbf{C}}}

\def\BF{{\mathbf{F}}}
\def\BG{{\mathbf{G}}}
\def\BH{{\mathbf{H}}}

\def\BL{{\mathbf{L}}}

\def\BP{{\mathbf{P}}}
\def\BQ{{\mathbf{Q}}}
\def\BR{{\mathbf{R}}}
\def\BS{{\mathbf{S}}}
\def\BT{{\mathbf{T}}}
\def\BU{{\mathbf{U}}}

\def\BZ{{\mathbf{Z}}}

\def\Bw{{\mathbf{w}}}

\def\Bpi{{\mathbf{\pi}}}

\def\eps{\varepsilon}

\def\Aut{\operatorname{Aut}\nolimits}

\def\can{{\mathrm{can}}}

\def\End{\operatorname{End}\nolimits}

\def\Ext{\operatorname{Ext}\nolimits}

\def\Gal{\operatorname{Gal}\nolimits}
\def\GL{\operatorname{GL}\nolimits}

\def\Hilb{\operatorname{Hilb}\nolimits}

\def\Ho{\operatorname{Ho}\nolimits}

\def\Hom{\operatorname{Hom}\nolimits}

\def\id{\operatorname{id}\nolimits}

\def\Ind{\operatorname{Ind}\nolimits}

\def\Irr{\operatorname{Irr}\nolimits}

\def\mMod{\operatorname{\!-mod}\nolimits}

\def\opp{{\operatorname{opp}\nolimits}}

\def\mperf{\operatorname{\!-perf}\nolimits}

\def\PSL{\operatorname{PSL}\nolimits}

\def\Out{\operatorname{Out}\nolimits}

\def\Res{\operatorname{Res}\nolimits}

\def\SL{\operatorname{SL}\nolimits}

\def\SU{\operatorname{SU}\nolimits}

\def\Spec{\operatorname{Spec}\nolimits}

\def\mstab{\operatorname{\!-stab}\nolimits}
\def\Stab{\operatorname{Stab}\nolimits}

\def\ie{{\em i.e.}}

\newtheorem{thm}{Theorem}[section]
\newtheorem{lemma}[thm]{Lemma}

\newtheorem{conj}[thm]{Conjecture}
\newtheorem{question}[thm]{Question}
\theoremstyle{definition}
\newtheorem{rem}[thm]{Remark}

\input xypic
\xyoption{all}
\xyoption{2cell}

\usepackage{tikz}
\usetikzlibrary{decorations.pathreplacing,calligraphy}

\title{Modular representations of finite groups and Lie theory}
\author{Rapha\"el Rouquier}
\address{Department of Mathematics, UCLA, Box 951555,
Los Angeles, CA 90095-1555, USA}
\email{rouquier@math.ucla.edu}
\dedicatory{In memory of Jacques Tits}
\thanks{The author gratefully acknowledges support from the NSF (grant DMS-1702305) and the
Simons Foundation (grant \#376202). This material is based upon work supported by a grant from the
Institute for Advanced Study.}

\begin{document}
\maketitle
\setcounter{tocdepth}{2}

\begin{abstract}
This article discusses the modular representation theory of finite groups of Lie type from the viewpoint
of Brou\'e's abelian defect group conjecture.
	We discuss both the defining characteristic case, the
inspiration for Alperin's weight conjecture, and the non-defining case, the inspiration for
Brou\'e's conjecture. The modular representation theory of general finite groups is conjectured to behave
both like that of finite groups of Lie type in defining characteristic, and in non-defining characteristic,
to a large extent.
	
	The expected behaviour of modular representation theory of finite groups of Lie
type in defining characteristic is particularly difficult to grasp along the lines of Brou\'e's conjecture,
	and we raise a new question 
related to the change of central character.

We introduce a degeneration method in the modular representation theory of finite groups
of Lie type in non-defining characteristic. Combined with the rigidity property of perverse equivalences,
this provides a setting for two-variable decomposition matrices, for large characteristic. This should
help make progress towards finding decomposition matrices, an outstanding problem with few general
results beyond the case of general linear groups. This last part is based on joint work with David Craven
and Olivier Dudas.
\end{abstract}

\tableofcontents

\section{Introduction}
Every finite simple group is a finite group of Lie type, an alternating group,
a cyclic group of prime order,  or one of the 26
sporadic groups. This provides a central role for finite groups of Lie type in
finite group theory.

Conjectures of Alperin and Brou\'e predict that the modular representation
theory of general finite groups shares many features with that of finite groups
of Lie type. Alperin's prediction is inspired by finite
groups of Lie type in defining characteristic. On the other hand, Brou\'e 
predicts a behaviour similar to that of finite groups of Lie type in
non-defining characteristic.

\medskip

For simple finite groups of Lie type in 
defining characteristic, the assumptions of Brou\'e's conjecture (abelian defect groups)
are only satisfied for groups of
type $A_1$, like $\PSL_2(\BF_q)$, outside cases of simple blocks. Brou\'e's conjecture is known to hold
in that
case, but the combinatorics involved in the proof have so far not
been understood within the usual Lie theoretic or geometric framework for
$\SL_2$.
An important problem is to find a proof of Brou\'e's conjecture for
$\SL_2(\BF_q)$ that relates to the geometry associated with the group. A
major open problem is to find an extension of Brou\'e's conjecture that 
removes the assumption on Sylow subgroups, and understanding Brou\'e's
conjecture for defining characteristic representations of $\SL_2(\BF_q)$
could lead to understanding how the conjecture should extend to higher
rank groups, and eventually to all finite groups.
Brou\'e's conjecture is about the existence of certain equivalences of
derived categories.

There are few known equivalences between blocks of finite groups of Lie type in defining characteristic. We could locate two: the derived equivalence between the principal block of $\mathrm{SL}_2(q)$ and the principal block of a Borel subgroup, and a similar result for the non-simple non-principal block when $p$ is odd.
In particular, the two non-simple blocks are derived equivalent. We propose to consider a generalization of this situation. 
A particular case of that extension applies to $G=\mathrm{SL}_r(q)$, $r$ a prime dividing $q-1$:
the non-simple blocks all have the same number of simple modules and the corresponding blocks for proper local subgroups are isomorphic. 

\medskip
Progress in the understanding of modular representations of finite groups
of Lie type in non-defining characteristic has mostly been achieved by
extending some of the work of Lusztig (and Deligne-Lusztig) about 
characteristic $0$ representations to characteristic $\ell$. Brou\'e's
conjecture in this case has a formulation in terms of Deligne-Lusztig theory.
The main difficulty in proving the conjecture  has been
about obtaining information about individual cohomology groups of
Deligne-Lusztig varieties, rather than about their alternating sum. In
particular, a key required vanishing property is still open in general,
even for characteristic $0$ coefficients.

\smallskip
The introduction of the notion of perverse equivalences and the conjecture that the equivalences expected
from Deligne-Lusztig varieties should be perverse (joint work with Joe Chuang)
lead to the fact that torus and Weyl group data determine the decomposition matrices for large enough
characteristic, using the conjectural combinatorial perversity function of Craven, and a rigidity
property of perverse equivalences (joint work with David Craven).

We explain two ways in which toroidal structures appear by degeneration in the modular representation
theory of finite groups of Lie type in non-defining characteristic.
We give a "global" topological construction using a limit of completed classifying spaces,
the starting point being Friedlander's description of the completed classifying space of a finite group of
Lie type in terms of
homotopy fixed points on the classifying space of the corresponding Lie group. We provide also
an explicit local algebraic construction. These lead to conjectural two-variable decomposition matrices
for large characteristic (joint work with Olivier Dudas).

\smallskip
In part \ref{se:finitegroups}, we consider general finite groups. We review $p$-local group
theory
and $p$-local representation theory and discuss Alperin and Brou\'e's conjectures. We introduce perverse
equivalences, a type of derived equivalences between derived categories with filtrations, that induces
abelian equivalences up to shifts on the slices of the filtration.

Part \ref{se:FiniteLie} introduces finite groups of Lie type as fixed points of a Frobenius endomorphism,
or a more
general Steinberg endomorphism of a reductive algebraic group. We discuss the $p$-local structure of
finite groups of Lie type, both in defining and non-defining characteristic.

Part \ref{se:defining} is devoted to the modular representation theory of finite groups of Lie type in
defining characteristic. We provide an explanation for Alperin's conjecture. In the case
of groups of finite Lie rank $1$, we discuss
the relation between representations of the group and of a Borel subgroup. We analyze next
when a block with a given central character has the same number of simple modules as the principal block
and raise the question of understanding the relation between the module categories of such blocks.

Parts \ref{se:nondef} and \ref{se:deg} are concerned with the modular representation theory of finite
groups of Lie type in non-defining characteristic $\ell$. We start with a discussion of Deligne-Lusztig 
varieties and endomorphisms coming from braid groups. We review Lusztig's theory of characteristic
$0$ representations and describe modular counterparts. We finish \S \ref{se:nondef} with a discussion
of the particular form of Brou\'e's conjecture for finite groups of Lie type in non-defining
characteristic.

Part \ref{se:deg} discusses two approaches to generic phenomena. The first approach is based on
the description of $\ell$-completed classifying spaces of finite groups of Lie type in terms of fixed
points under unstable Adams operations and we discuss the rigidification of a certain limit of those
Adams operation, in relation with classifying spaces of loop groups. The second approach is based
on a degeneration of the group algebra of the local block, and the relation with the rigidity of perverse
equivalences. For general linear and unitary groups, there is a further relation with Hilbert schemes
of points on surfaces.

We gather in the appendix a number of basic facts on representations of algebras and finite groups,
in particular in relation with various types of equivalences. We give a very succinct survey of
basic constructions involving complex reflection groups, braid groups and Hecke algebras.

\medskip
I thank C\'edric Bonnaf\'e, David Craven, Olivier Dudas, Jesper Grodal, George Lusztig and Gunter Malle
for their comments.

\smallskip
This article is based on two lectures given at the Institute for Advanced Study,
Princeton, in the fall 2020.

\section{Finite groups}
\label{se:finitegroups}

\subsection{Group theory}
\subsubsection{Classification of finite simple groups}

Every finite simple group is one of the following \cite[\S 47]{Asch}
\begin{itemize}
	\item a cyclic group of prime order
	\item an alternating group $\GA_n$ for $n\ge 5$
	\item a finite simple group of Lie type
	\item one of the 26 sporadic simple groups, with orders ranging from $7920$ for the Mathieu group
		$M_{11}$ discovered in 1861 to about $8\times 10^{53}$ for the Fischer-Griess monster
		discovered in 1973.
\end{itemize}
Finite groups of Lie type govern to a large extent the structure of general finite groups.

\subsubsection{$p$-local group theory}
Let $G$ be a finite group, $p$ a prime and $P$ a Sylow $p$-subgroup of $G$.

\smallskip
The $p$-local group theory is the study of $G$ using {\em $p$-local subgroups}, i.e. subgroups of the form
$N_G(Q)$ for $Q$ a non-trivial $p$-subgroup.

This was originally developed mostly in the case $p=2$ and this 
underlies the proof of the classification of finite simple groups.

\smallskip
Here are some examples of results of $p$-local group theory. We denote by $O_{p'}(G)$ the largest
normal subgroup of $G$ of order prime to $p$.

\begin{itemize}
\item (Burnside 1897) When $P$ is abelian, two elements of $P$ that are conjugate in $G$ are also
		conjugate in $N_G(P)$ \cite[Chap.7, Theorem 1.1]{Go}.
\item (Frobenius) If $N_G(Q)/C_G(Q)$ is a $p$-group for every non trivial subgroup $Q$ of $P$, then 
	$G=O_{p'}(G)\rtimes P$ \cite[Chap. 7, Theorem 4.5]{Go}.
\item (Brauer-Fowler 1955) If $G$ is simple and $s\in G$ is an involution, then $|G|\le (2|C_G(s)|^2)!$ \cite[(45.4)]{Asch}.
\item (Brauer-Suzuki 1959) If the Sylow $2$-subgroups of $G$ are quaternion groups, then $G$ is not
	simple \cite[\S 3.3]{Co}.
\item (Glauberman 1966, case $p=2$) If $x$ is an element of order $p$ of $P$ that is not $G$-conjugate to
	any other element of $P$, then $G=O_{p'}(G)C_G(x)$ \cite[Appendix]{Co}.
\item (Alperin 1967) One can tell if two elements of $P$ are conjugate in $G$ using only $p$-local subgroups
	\cite[(38.1)]{Asch}.
\end{itemize}

Glauberman's Theorem (which generalizes Brauer-Suzuki's Theorem) 
holds also for odd primes, but the proof for odd primes uses the classification of
finite simple groups. Modular representation theory of finite groups was developed by Brauer as a tool
for studying finite groups. For example, the proof of the Brauer-Suzuki Theorem uses representation
theory in characteristic $2$.

It is hoped that modular representation theory will eventually reach a point where it can be used to
obtain a direct proof of Glauberman's Theorem for odd primes and lead to simplifications of the proof
of the classification of finite simple groups.

Modular representation theory leads to a generalization of local group theory, where Sylow subgroups
are replaced by defect groups of blocks. A major theme of modular representation theory is to relate
modular representations of a group and its local subgroups, sometimes leading to versions of
"factorization" results like Frobenius's Theorem above replaced by an equivalence between module
categories \cite{BrouPu}.

\subsection{$p$-local representation theory}

\subsubsection{}
$p$-local representation theory is concerned with the
study of representations of $G$ over $\BZ_p$ and $\BF_p$ (or finite extensions of those) in relation with
$p$-local subgroups and their representations.

It involves character-theoretic aspects, in particular the value of
complex characters of $G$ on elements whose order is divisible by $p$. It involves also the study of
simple and indecomposable representations and mod-$p$ cohomology.

\subsubsection{}
Let $p$ be a prime number and $\CO$ be the ring of integers of a finite extension $K$ of $\BQ_p$. Let $k$
be the residue field of $\CO$.

Let $G$ be a finite group. The category $kG\mMod$ is not semisimple if $p$ divides $|G|$, but it still
splits as direct sum of indecomposable
full abelian subcategories. This is induced by a corresponding decomposition
of $\CO G\mMod$. That decomposition comes from a decomposition of $1$ as a sum of orthogonal
primitive idempotents
$1=\sum_b b$ of $Z(\CO G)$, the block idempotents. We have an algebra decomposition into blocks
$\CO G=\prod_b \CO Gb$ and a category decomposition $\CO G\mMod=\bigoplus_b \CO Gb\mMod$.
We will still denote by $b$ the image of the idempotent in $Z(kGb)$ and we have corresponding
decompositions of $kGb$ and $kGb\mMod$.

The {\em principal block} $\CO Gb_0$ of $\CO G$ is the one such that $b_0$ does not act by $0$ on the
trivial representation.

\smallskip
We will always assume that $K$ contains all $|G|$-th roots of unity.

\medskip
A {\em defect group} of a block $\CO Gb$ is a minimal subgroup $D$ of $G$ such that the restriction functor
$D^b(\CO Gb)\to D^b(\CO D)$ is faithful. A defect group is a $p$-subgroup of $G$ and all
defect groups are conjugate. The defect groups are trivial
if and only if $\CO Gb$ is a matrix algebra over $\CO$ (equivalently, $kGb$ is semisimple). The
defect groups of the principal block are the Sylow $p$-subgroups of $G$.

There is a unique block idempotent $b_D$ of $\CO N_G(D)$, the {\em Brauer correspondent} of $b$,
such that the functor 
$b_D\CO Gb\otimes_{\CO Gb}-:D^b(\CO Gb)\to D^b(\CO N_G(D)b_D)$ is faithful. The idempotent
$b_D$ is actually contained in $\CO C_G(D)$.

\subsection{Conjectures}
\subsubsection{Alperin's weight conjecture}

Alperin's weight conjecture \cite{Alp} asserts that the
number of non-projective simple $kG$-modules is locally determined.

A {\em weight} for $G$ is a pair $(Q,V)$, where $Q$ is a $p$-subgroup of $G$ and
$V$ is a projective simple $kN_G(Q)/Q$-module.

\begin{conj}[Alperin]
	The number of isomorphism classes of simple $kG$-modules is the same as the number of
conjugacy classes of weights.
\end{conj}

The conjecture has also a blockwise version.

\smallskip
We will explain in \S \ref{se:AlperinLie}
that the conjecture has a bijective proof when $G$ is a finite group of Lie type
in defining characteristic, for example $G=\GL_n(\BF_{p^r})$. This is the inspiration for the conjecture.
In a sense, Alperin's conjecture predicts
that all finite groups behave like finite group of Lie type in defining characteristic.

\subsubsection{Brou\'e's conjecture}
Brou\'e's conjecture \cite{Brou1} asserts that the derived category of
$kG$-modules (excluding as above semi-simple parts) is determined locally, when Sylow $p$-subgroups are
abelian. There is a blockwise version which we now state.

\begin{conj}[Brou\'e]
	\label{co:Broue}
	Let $\CO Gb$ be a block with abelian defect group $D$. We have
	$D^b(\CO Gb)\simeq D^b(\CO N_G(D)b_D)$.
\end{conj}

When $G$ is a finite group of Lie type in non-defining characteristic, Brou\'e and others have proposed
an explicit candidate for a functor realizing an equivalence, using Deligne-Lusztig varieties
(cf \S\ref{se:BroueconjP}).
In a sense, Brou\'e's conjecture predicts
that all finite groups behave like finite groups of Lie type in non-defining characteristic as far as
modular representations are concerned (with the abelian defect assumption), even though
the sought-after equivalence won't arise from something like a Deligne-Lusztig variety.

\begin{rem}
Brou\'e's conjecture does not extend in an obvious way to blocks with non-abelian defect groups, cf for
example \S\ref{se:rank1}.

On the other hand, one can generalize slightly Brou\'e's conjecture to the case where the hyperfocal
	subgroup of the defect group is abelian \cite[Appendix A.2]{Rou1}.
\end{rem}

\subsubsection{Comparison}
The two conjectures lead to an odd phenomenon: finite groups, with respect to a prime $p$, behave
like finite groups of
Lie type in both defining {\em and} non-defining characteristic!

\smallskip
For blocks with abelian defect groups, Brou\'e's conjecture implies Alperin's conjecture.
A major open problem (beyond proving these conjectures) is to find a structural statement like Brou\'e's
conjecture for general blocks.

\smallskip
Note that the neighborhood of the trivial representation is determined locally: $H^*(G,k)$ can be
recovered as a subalgebra of $H^*(P,k)$.
When $P$ is abelian (or when $N_G(P)$ controls fusion), then 
$H^*(G,k)=H^*(P,k)^{N_G(P)}=H^*(N_G(P),k)$.

A general version of Brou\'e's conjecture should 
contain both that fact and the information about the number of simple modules.
Going beyond the neighborhood of the trivial representation is discussed for finite groups
of Lie type in non-defining characteristic in \S\ref{se:rigidification}.

\medskip
Alperin's conjecture can be reformulated in terms of chains of $p$-subgroups \cite{KnRob}.
Consider the poset of non-trivial $p$-subgroups of $G$ and the associated simplicial complex $\CP$
($n$-simplices are chains $Q_0<Q_1<\cdots<Q_n$).

Alperin's conjecture (for all finite groups) is equivalent to 
the equality (for all finite groups $G$)
$$\sum_{c\in\CP/G} (-1)^{|c|}l'(k\Stab_G(c))=0$$
where $l'(kG)$ is the number of non-projective simple modules.
The simplicial complex $\CP$ can be replaced by its subcomplex given by elementary abelian subgroups or by
other complexes using subgroups $Q$ such that $Q=O_p(N_G(Q))$: the sum does not change. Here
we denote by $O_p(H)$ the largest normal $p$-subgroup of a finite group $H$.
 
There is also a blockwise version of that reformulation.
This reformulation suggests that appropriate categories of representations of local subgroups could
be glued to recover some weakened version of the category of representations of $G$.

\subsection{Perverse equivalences}

\subsubsection{Definition}
Let $A$ and $A'$ be two finite dimensional algebras over a field $k$ and let
$F:D^b(A)\iso D^b(A')$ be an equivalence of triangulated categories.
The equivalence $F$ induces an isomorphism of abelian groups
$K_0(D^b(A))=\BZ^{\Irr(A)}\iso \BZ^{\Irr(A')}=K_0(D^b(A'))$, but no bijection $\Irr(A)\iso\Irr(A')$.
So, in the situation of Brou\'e's abelian defect conjecture (Conjecture \ref{co:Broue}), there is
no expectation of a bijection between simple modules. We will introduce now a particular
type of derived equivalence that induces such a bijection \cite{ChRou}.

\smallskip
Fix $\pi:\Irr(A)\to\BZ$. 
An equivalence $F:D^b(A)\iso D^b(A')$ is {\em perverse} relative to $\pi$ if
there is a bijection $f:\Irr(A)\iso \Irr(A')$ such that
\begin{itemize}
	\item given $S\in\Irr(A)$, if $T$ is a composition factor of $H^i(F(S))$,
		then $\pi(f^{-1}(T))<\pi(S)$ for $i\neq -\pi(S)$
	\item $H^{-\pi(S)}(F(S))$ admits $f(S)$ as a composition factor with multiplicity one, and all
		other composition factors $T$ satisfy $\pi(f^{-1}(T))<\pi(S)$.
\end{itemize}
When this holds, the map $f$ is determined by $F$ and $\pi$.

\smallskip
Given $A$ and $\pi$, then $A'$ is unique up to Morita equivalence (if it exists).

\begin{rem}
This is a particular case of a more general definition that involves the additional data of
	an order on $\Irr(A)$.

Perverse equivalences can be defined more generally for derived categories of abelian categories. A further
generalization is the consideration of a filtered triangulated category with two $t$-structures and the
notion of a shift of $t$-structures with respect to a perversity function: the $t$-structures
are assumed to be compatible with the filtration and the $t$-structures induced on the slices
of the filtration differ by a shift given by the perversity function.
\end{rem}

\subsubsection{Examples}
\label{se:pervexamples}

Consider the situation of Brou\'e's conjecture: we have a block $\CO Gb$ with defect group $D$.
If there is an equivalence $D^b(\CO Gb)\iso D^b(\CO N_G(D)b_D)$ such that the induced equivalence over
$k$ is perverse, then there is a total order on $\Irr(kGb)$ and on $\Irr(KGb)$ such that 
the decomposition matrix of $\CO Gb$ has the following shape (cf \S\ref{se:appfinite} for the definition of
decomposition matrices):
$$\begin{matrix}
1 \\
*&1 \\
\vdots&\ddots&\ddots \\
*&\hdots&*& 1 \\
	*&\hdots&\hdots&* \\
	\vdots&&&\vdots \\
	*&\hdots&\hdots&* 
\end{matrix}
$$

\smallskip
In the situation of Brou\'e's abelian defect conjecture, 
perverse equivalences are known to exist in a number of cases: when defect groups are cyclic 
and conjecturally in the case of groups of Lie type in non-defining characteristic \cite{ChRou}
(cf \S \ref{se:BroueconjP}), for principal blocks with
defect group of order $4$ or $9$ of groups with no simple factor $\GA_6$ or $M_{22}$ when $p=3$ 
\cite[Theorem 4.36]{CrRou1}.
Note that there is no perverse equivalence in the situation of Brou\'e's conjecture for 
principal $3$-blocks of $\GA_6$ and $M_{22}$ \cite[\S 5.3.2.3 and \S 5.4.3]{CrRou1}.

\section{Finite groups of Lie type}
\label{se:FiniteLie}

We discuss now finite groups of Lie type (cf \cite{MalTe}).

\subsection{Reductive groups}
\label{se:reductive}

Let $p$ be a prime number and $\bar{\BF}_p$ be an algebraic closure of
the finite field with $p$ elements. Given $q$ a power of $p$, we denote by $\BF_q$ the subfield of
$\bar{\BF}_p$ with $q$ elements.

\smallskip
Let $\BG$ be a (connected) reductive (linear) algebraic group over $\bar{\BF}_p$.

Let $\BT_0$ be a maximal torus of $\BG$ and $\BB_0$ a Borel subgroup of $\BG$ containing $\BT_0$.
Let $X=X(\BT_0)=\Hom(\BT_0,\BG_m)$ and $Y=Y(\BT_0)=\Hom(\BG_m,\BT_0)$.
Let $\Phi\subset X$ denote the set of roots,
$\Delta$ the set of simple roots, $\Phi^\vee$ the set of coroots and
$\Delta^\vee$ the set of simple coroots. Given
$\alpha\in\Delta^\vee$, we denote by $\omega_\alpha$ the corresponding fundamental weight.
Let $\rho$ be the half sum of the positive roots.

 Let $W=N_{\BG}(\BT_0)/\BT_0$ be
the Weyl group and let $S=(s_\alpha)_{\alpha\in\Delta}$ be its generating set as a Coxeter group.

\subsection{Rational structures}
\subsubsection{Frobenius endomorphisms}
Let $V_0=\Spec A_0$ be an affine algebraic variety over $\BF_q$, where $q$ is a power of $p$.
The endomorphism $F$ of $V=V_0\times_{\BF_q}\bar{\BF}_q=\Spec (A_0\otimes_{\BF_q}\bar{\BF}_p)$
given on $A_0\otimes_{\BF_q}\bar{\BF}_p$ by $a\mapsto a^q$ is called the (geometric) Frobenius endomorphism
of $V$.

Given $V'$ an affine algebraic variety over $\bar{\BF}_p$, an endomorphism $F$ of $V'$ is
called a {\em Frobenius endomorphism} if there is a power $q$ of $p$ and
an affine algebraic variety $V_0$ defined over $\BF_q$ with an isomorphism 
$V'\iso V_0\times_{\BF_q}\bar{\BF}_q$ identifying $F$ with the Frobenius endomorphism coming from $V_0$.
We say that $F$ defines a {\em rational structure} for $V'$ over $\BF_q$.

\subsubsection{Steinberg endomorphisms}
\label{se:Steinberg}

Let $F$ be an endomorphism of $\BG$, a power of which is a Frobenius endomorphism (such an $F$ is called
a {\em Steinberg endomorphism} of $\BG$).
The group $G=\BG^F$ is finite: this is a finite group of Lie type.
We will put more generally $H=\BH^F$ for $\BH$ an $F$-stable subgroup of $\BG$.

There exists
an $F$-stable Borel subgroup $\BB_0$ and an $F$-stable maximal torus $\BT_0$ contained in $\BB_0$.

\smallskip

We say that $(\BG,F)$ is {\em split} if $F$ acts by a multiple of the identity on
$X(\BT_0)$. If $(\BG,F)$ is split then $F$ is a Frobenius endomorphism.

\smallskip
Let $\delta$ be the minimal positive integer such that $(\BG,F^\delta)$ is split and we 
define the positive real number $q$ by the requirement that $F^\delta$  defines
a rational structure over $\BF_{q^\delta}$.
The automorphism 
of $X\otimes_\BZ\BQ$ induced by $F$ permutes the lines $\BQ\alpha$ for $\alpha\in\Delta$ and
this provides an automorphism $\phi$ of order $\delta$ of the Coxeter diagram of $\BG$. When
$\BG$ is simple and simply connected, the
endomorphism $F$ of $\BG$ depends only on $q$ and on $\phi$, up to an inner automorphism.

If $F$ is a Frobenius endomorphism, then $F(\alpha)=q\phi(\alpha)$ for all $\alpha\in\Delta$.

\smallskip
Let $\bar{\delta}$ be the minimal positive integer such that $F^{\bar{\delta}}$ is a Frobenius
endomorphism of $\BG$.

\subsection{Finite simple groups of Lie type}
The finite simple groups of Lie type are obtained by the following construction. We assume that $\BG$ is
simple and simply connected. Then, the group $G/Z(G)=\BG^F/Z(\BG)^F$ is simple (with some exceptions described after
the classification below) and we obtain in this way all
finite simple groups of Lie type. The group is denoted by ${^\delta D}(q^{\bar{\delta}})$, where $D$ is
the type
of the root system of $\BG$ (equivalently, the type of the Dynkin diagram). We give now the list of all
those finite simple groups, with restrictions on $q$ if any.

\begin{itemize}
	\item $A_n(q)$, $n\ge 1$
	\item $B_n(q)$, $n\ge 2$
	\item $C_n(q)$, $n\ge 3$
	\item $D_n(q)$, $n\ge 4$
	\item $E_6(q)$, $E_7(q)$, $E_8(q)$, $F_4(q)$, $G_2(q)$
	\item ${^2A}_n(q)$, $n\ge 2$
	\item ${^2D}_n(q)$, $n\ge 3$
	\item ${^3D}_4(q)$, ${^2E}_6(q)$
	\item ${^2B}_2(q^2)$, $q^2=2^{2m+1}$ for some $m\ge 1$
	\item ${^2G}_2(q^2)$, $q^2=3^{2m+1}$ for some $m\ge 1$
	\item ${^2F}_4(q^2)$, $q^2=2^{2m+1}$ for some $m\ge 0$
\end{itemize}
with the following exceptions: $A_1(2)$, $A_1(3)$, ${^2A}_2(2)$ and  $B_2(2)$ are not simple
and can be removed from the list.
The derived subgroup of ${^2F}_4(2)$ is simple and of index $2$ in ${^2F}_4(2)$. It does not arise in another construction and that group
(the Tits group) needs to be added to the list of finite simple groups of Lie type.
Note that some finite simple groups of Lie type occur more than once in the classification.

\subsection{Local structure}
We consider the setting of \S \ref{se:Steinberg}.

\subsubsection{Defining characteristic}
\label{se:localp}

Let $\BU_0$ be the unipotent radical of $\BB_0$. Then $U_0$ is a Sylow $p$-subgroup of $G$.

An important class of subgroups in local group theory consists of those $p$-subgroups $Q$ of $G$ such that
$Q=O_p(N_G(Q))$. In our setting where $G$ is a finite group of Lie type
and $p$ is the defining characteristic, $Q$ satisfies this
condition if and only if $Q$ is the group of $F$-fixed points of the unipotent radical of an
$F$-stable parabolic subgroup of $\BG$ \cite[Remark 6.15]{CaEn2}.

\medskip
Note that $U_0$ is abelian if and only if the quotient of $\BG$ by its radical is a product of groups
of type $A_1$.

\begin{rem}
Assume $\BG$ is simple and simply connected 
with $F$-rank $1$, \ie\ of type $A_1$ (group $\SL_2(q)$), ${^2A}_2$ (group
$SU_3(q^2)$), ${^2B}_2$ (Suzuki group) or ${^2 G}_2$ (Ree group).

In that case, the Sylow $p$-subgroups of $G$ have the trivial intersection property: given $g\in G$ with
$g{\not\in}N_G(U_0)=B_0$, we have $U_0\cap gU_0g^{-1}=1$. This implies that two subgroups of $U_0$ that are
conjugate in $G$ are already conjugate in $B_0$.
\end{rem}

\subsubsection{Non-defining characteristic}
\label{se:localnondef}
We consider a prime number $\ell\neq p$.

We denote by $d$ the order of $q$ in $\BF_\ell^\times$. 

\smallskip
We assume here that $F$ is a Frobenius endomorphism. The $\ell$-local structure of $G$ has
generic aspects explained in \cite{BrouMal}.

\smallskip
A {\em $\Phi_d$-subgroup} of $\BG$ is an $F$-stable torus $\BS$ such that $\Phi_d(q^{-1}F)$ acts by $0$ on
$Y(\BS)\otimes\BQ$.  Here, $\Phi_d$ is the $d$-th cyclotomic polynomial.

\smallskip
Let $\BS$ be a $\Phi_d$-subgroup of $\BG$.
There is $g\in\BG$ such that $g\BT_0g^{-1}$ is $F$-stable and contains $\BS$. 
Let $w$ be the image of $g^{-1}F(g)\in N_{\BG}(\BT_0)$ in $W$.
Given $s\in\BS$, we have
$wF(g^{-1}sg)w^{-1}=g^{-1}F(s)g$. It follows that the image of $Y(g^{-1}\BS g)\otimes\BQ$ in
$Y(\BT_0)\otimes\BQ$ is a subspace
on which $\Phi_d(w(q^{-1}F))$ acts by $0$.

All maximal $\Phi_d$-subgroups of $\BG$ are $G$-conjugate. 
Let $\BS$ be a maximal $\Phi_d$-subgroup of $\BG$ and $g$, $w$ be as above.
The image of $Y(\BS)$ in $Y(\BT_0)$ is $Y(\BT_0)\cap\ker(\Phi_d(wq^{-1}F))\ |\ Y(\BT_0)\otimes\BQ)$.

\medskip
Assume $\ell>3$ and $F$ is a Frobenius endomorphism. Let
$\BS$ be a maximal $\Phi_d$-subgroup of $\BG$ and $\BL=C_\BG(\BS)$, a Levi subgroup.
Let $\BT_0$ be a maximal $F$-stable
torus of $\BL$ and let $W_\BL=N_{\BL}(\BT_0)/\BT_0$.
Let $W'$ be the subgroup of $W$
generated by roots orthogonal to all the roots corresponding to $\BL$.
There exists an $\ell$-subgroup $D'$ of $N_\BG(\BT_0)^F$ whose image in $W$ is a Sylow $\ell$-subgroup
of $W'$ and such that $D'\cap T_0=1$.
Then $D=Z(\BL)^F_\ell\rtimes D'$ is a Sylow $\ell$-subgroup of $G$ \cite{Ca2} (cf also
\cite[Exercice 22.6]{CaEn2}).
Furthermore, $D$ is abelian if and only if $D'=1$.

\smallskip
When $D$ is abelian, we have an isomorphism $N_G(D)/C_G(D)\iso (N_W(W_\BL)/W_\BL)^{wF}$, where $w$ is
defined as above from $\BS$. The group $N_G(D)/C_G(D)$ is a complex reflection group
\cite[Theorem 3.4]{LeSp} (cf \cite[Theorem 5.7]{Brou3}).
When $\BL$ is a torus, then $q^{-1}wF$ is $d$-regular and this is explained in \S\ref{se:regular}.

\begin{rem}
	It is a remarkable fact that when $D$ is abelian, then $N_G(D)/C_G(D)$ is a reflection group.
	We showed
	in \cite{Rouauto} that a suitable version of this property actually holds for all finite simple
	groups.
	We consider $G_0$ a simple group with an abelian Sylow $\ell$-subgroup $D$ such that 
	$H^2(G,\BF_\ell)=0$. Let $G$ be a subgroup of $\Aut(G_0)$ containing $G_0$ and
	such that $G/G_0$ is a Hall $\ell'$-subgroup of $\Out(G_0)$.  Then there exists 
	\begin{itemize}
		\item a field extension $K$ of $\BF_\ell$,
		\item an extension of the structure of $\BF_\ell$-vector space on the largest elementary
			abelian subgroup $\Omega_1(D)$ of $D$ to a structure of $K$-vector space
		\item a subgroup $N$ of $\GL_K(\Omega_1(D))$
		\item a subgroup $\Gamma$ of $\Aut(K)$
	\end{itemize}
	such that $N_G(D)/C_G(D)=N\rtimes\Gamma$ and the normal subgroup of $N$ generated by reflections
	acts irreducibly on $\Omega_1(D)$.

	For example, when $G_0=\PSL_2(\ell^n)$, we view $D\simeq (\BZ/\ell)^n$ as a one-dimensional
	vector space over $K=\BF_{\ell^n}$, we have $N=K^\times$ and
	$\Gamma$ is a Hall $2'$-subgroup of $\Gal(K/\BF_\ell)$.
	
	\smallskip
	It would be very interesting to find a role for reflection groups in Brou\'e's abelian
	defect group conjecture (Conjecture \ref{co:Broue}).
\end{rem}

\section{Defining characteristic}
\label{se:defining}
As we will see, it is much easier to parametrize irreducible representations of $G$ in characteristic $p$
than in characteristic $0$ (cf \cite{Hu}).

\subsection{Simple modules and blocks}
\subsubsection{Rational representations}

Let $\BG$ be a reductive connected algebraic group over an algebraic closure $\bar{\BF}_p$ of
a finite field with $p$ elements, where $p$ is a prime number. We consider $\BT_0$, $\BB_0$, etc as in
\S\ref{se:reductive}.

Let $X_+(\BT_0)=\{\lambda\in X(\BT_0)\ |\ \langle \lambda,\alpha^\vee\rangle\ge 0\ \forall\alpha\in\Delta\}$
be the set of dominant weights.
Given $\lambda\in X_+(\BT_0)$, let $\CL_\lambda=\BG\times_{\BB_0}(\bar{\BF}_p)_\lambda$ be the associated line bundle
on the flag variety $\BG/\BB_0$.
The rational $\BG$-module $H^0(\BG/\BB_0,\CL_\lambda)$ has a unique simple submodule
$\BL(\lambda)$ and $\{\BL(\lambda)\}_{\lambda\in X_+(\BT_0)}$ is a full set of representatives
of isomorphism classes of simple rational $\BG$-modules.

\subsubsection{Representations of the finite group}
\label{se:repatp}
We consider now the setting of \S \ref{se:Steinberg}.

Given $r\ge 0$, we put $X_r=\{\lambda\in X_+(\BT_0)\ |\ \langle\lambda,\alpha^\vee\rangle<r\ \forall\alpha\in
\Delta\}$.

\smallskip

Given $\lambda\in X_+$, let $L(\lambda)$ be the restriction of $\BL(\lambda)$ to
$G$.

Given $A$ an abelian group, we put $A^\vee=\Hom(A,\bar{\BF}_p^\times)$.

\smallskip

We assume in \S\ref{se:repatp} that $\BG$ is simply connected and simple and that
$F$ is a Frobenius endomorphism.

\smallskip
We have the following description of simple modules (Steinberg) \cite[Theorem 2.11 and \S 20.2]{Hu}
and of their blocks and defect groups (Dagger and Humphreys) \cite[Theorem 8.5 and \S 20.3]{Hu}.

\begin{thm}
	\label{th:simplesdescribing}
	The set $\{L(\lambda)\}_{\lambda\in X_q}$ is a complete set of representatives
	of isomorphism classes of simple $kG$-modules.

	There is one block of defect zero, with simple module the Steinberg module 
	$L((q-1)\rho)$. The other blocks have maximal defect, they are parametrized by 
	$Z(G)^\vee$. The simple modules in the block $kGb_\zeta$ corresponding to
	$\zeta\in Z(G)^\vee$ are the $L(\lambda)$ with $\lambda_{|Z(G)}=\zeta$ 
	and $\lambda\in X_q\setminus\{(q-1)\rho\}$.
\end{thm}

\begin{rem}
	Note that the set of simple $kG$-modules and their dimensions depend only on $q$, not on the
	Frobenius endomorphism. For example, the sets of simple modules for
	$\mathrm{SU}_n(q)$ and $\mathrm{SL}_n(q)$ are obtained by restricting the same set of
	simple rational representations of $\SL_n$.
\end{rem}

\begin{rem}
	When $F$ is not a Frobenius endomorphism, Theorem \ref{th:simplesdescribing} needs to be modified
	as follows \cite[\S 20]{Hu}.
	We assume $F$ is not a Frobenius endomorphism. Note that $F^2$ is a Frobenius
	endomorphism defining a rational structure over $\BF_{q^2}$. Note also that
	$Z(G)=1$. Define now $X_{q^2,S}$ as the set of $\lambda\in X_{q^2}$ such that 
	$\langle\lambda,\alpha^\vee\rangle=0$ for every long simple root $\alpha$.
	The set $\{L(\lambda)\}_{\lambda\in X_{q^2,S}}$ is a complete set of representatives
	of isomorphism classes of simple $kG$-modules.
	
	The Steinberg module is $L((q^2-1)\rho_S)$, where
	$$\rho_S=\sum_{\substack{\alpha\in S \\ \alpha\textrm{ short}}}\omega_\alpha.$$
	It is in a block of defect zero and the
	principal block is the unique other block.
\end{rem}

\subsection{Alperin's conjecture}

\subsubsection{Bijective proof}
\label{se:AlperinLie}
We assume in \S\ref{se:AlperinLie} that $\BG$ is 
simply connected and $F$ is a Frobenius endomorphism. We follow \cite{Ca1} and \cite[\S 6.3]{CaEn2}.

Given $I$ an $F$-stable subset of $S$, let $\BL_I$ be the corresponding standard Levi subgroup of $\BG$ and
let $X_I$ be the set of $\lambda\in X_q$ such that given $\alpha\in\Delta$, we have
$\langle\lambda,\alpha^\vee\rangle=q-1$ if and only if
$\alpha\in I$. We have $X_q=\coprod_{I\subset S}X_I$.

\medskip
$\bullet\ $Restriction from $\BT_0$ to $T_0$ induces a bijection $X_I\iso 
\Irr_k(T_0/(T_0\cap [L_I,L_I]))\iso L_I^\vee$.

\smallskip
$\bullet\ $Given $\zeta\in L_I^\vee$, the $kL_I$-module $\mathrm{St}_{L_I}\otimes\zeta$ is 
simple and projective. This provides a bijection from
$L_I^\vee$ to the set of isomorphism classes of projective simple $kL_I$-modules.

\smallskip
$\bullet\ $Let $Q$ be a $p$-subgroup of $G$ such that $kN_G(Q)/Q$ has a simple projective module. Then
$O_p(N_G(Q)/Q)=1$ and it follows (cf \S \ref{se:localp}) that
$Q$ is the subgroup of $F$-fixed points
of the unipotent radical of an $F$-stable parabolic subgroup of $\BG$. So, we have
a bijection from the union over $I$ an $F$-stable subset of $S$ of the sets of 
isomorphism classes of projective simple $kL_I$-modules to the set of
$G$-conjugacy classes of pairs $(Q,V)$ where $Q$ is a $p$-subgroup of $G$ and $V$ a simple
projective $kN_G(Q)/Q$-module, taken up to isomorphism.

\smallskip
Together with the bijection from $X_q$ to $\Irr_k(G)$, we obtain a bijection between
$\Irr_k(G)$ and  the set of
$G$-conjugacy classes of pairs $(Q,V)$ where $Q$ is a $p$-subgroup of $G$ and $V$ a simple
projective $kN_G(Q)/Q$-module, taken up to isomorphism. This confirms Alperin's conjecture for $G$.

$$\xymatrix{
	X^+\ar[dd]_\sim & X_q^+\ar[dd]_\sim\ar[l] \ar[r]^-\sim & \displaystyle \coprod_{\substack{
		I\subset S\\ I\ F\text{-stable}}}\Irr_k(T_0/(T_0\cap [L_I,L_I]))\ar[d]^\sim \\
&& \displaystyle \coprod_{\substack{
		I\subset S\\ I\ F\text{-stable}}}\{\text{proj simple }kL_I\text{-modules}\}\ar[d]^\sim \\
\Irr(\BG) &  \Irr_k(G) \ar[l] \ar@{-->}[r]_-\sim & 
\{(Q,V)\ | Q\le G,\ V \text{ proj simple }k(N_G(Q)/Q)\text{-module}\}/G\\
}$$

\subsubsection{Abelian defect}
\label{se:abelianp}
If $\BG$ is simple and $kG$ has a block with non-trivial abelian defect groups,
then $\BG$ is of type $A_1$ (cf \S \ref{se:localp} and Theorem \ref{th:simplesdescribing}).
Brou\'e's abelian defect group conjecture has been solved for $\SL_2(\BF_q)$
(\cite{Ch} for $q=p^2$, \cite{Ok} for the principal block, \cite{Yo} for the non-principal block with
maximal defect and \cite{Wo2} for the proof that the equivalence is a composition of perverse equivalences).
The solution involves some rather complicated combinatorial and algebraic constructions.

\begin{rem}
Assume $\BG$ is semisimple and simply connected, split over $\BF_p$ and assume that $p$ is larger
than the Coxeter number of $\BG$. We denote by $R=\bigoplus_{\alpha\in\Delta}\BZ\alpha$ the root
lattice.

Consider the full subcategory $\CC$ of the derived category of bounded complexes of
	finite-dimensional $\BB_0$-modules whose objects are those $C$ such that $H^i(C)$ has
weights in $pR$ for all $i$.
The functor $R\Ind_{\BB_0}^{\BG}$ induces an equivalence from $\CC$ to
the bounded derived category of the principal block of finite-dimensional representations of $\BG$
(\cite{ArkBeGi} for the case of quantum groups at a root of unity and 
\cite{HoKaSc} for an adaption and details of the characteristic $p$ case).

It would be very interesting if this equivalence could be used to relate representations
of $G$ and $B_0$ over $k$, particularly in the case of $\SL_2$. This would possibly shed light on how
Brou\'e's conjecture could be generalized to non-abelian defect groups.
\end{rem}

\subsubsection{Groups of $F$-rank $1$}
\label{se:rank1}
When $\BG$ has $F$-rank $1$, the induction and restriction functors provide inverse stable
equivalences (cf \S\ref{se:equivalences}) between $kB_0$ and $kG$ because the Sylow $p$-subgroups of $G$ have the trivial intersection
property (cf \S\ref{se:localp}).

The principal blocks of $kB_0$ and $kG$ are actually derived equivalent when $G=\SL_2(q)$, but this
is known not to generalize to all groups of $F$-rank $1$.
It is known in a number of cases that the principal blocks of $kG$ and $kB_0$ are not derived equivalent
because the centers are not isomorphic: for ${^2G}_2(q^2)$ when $q^2\ge 27$ \cite{BrougSchw} and for
$G=\SU_3(q)$ when $3\le q\le 8$ \cite{BouZi}. It is also known that the principal blocks of
$\CO G$ and $\CO B_0$ are not derived equivalent because the centers are not isomorphic
for $G={^2B}_2(q^2)$ when $q^2\ge 8$ \cite{Cl}. Not that in this last case the centers of the principal
blocks of $kG$ and $kB_0$ are isomorphic. It is expected that the principal blocks of $kG$ and $kB_0$ are
not derived equivalent in that case.

\subsection{Change of central character}
\label{se:change}
\subsubsection{Number of simple modules in a block}
We assume in \S\ref{se:change} that $F$ is a Frobenius endomorphism.

Let $\bar{\gamma}=\Res^{\BT_0}_{Z(G)}:X\to Z(G)^\vee$.
\begin{lemma}
	\label{le:numbersimple}
	Let $\zeta\in Z(G)^\vee$. The number of simple modules in $kGb_\zeta$
	is less than or equal to that for the principal block $kGb_0$.
	There is equality if and only if 
	$\zeta\in \bigcap_{\alpha\in\Delta^\vee} \bar{\gamma}(\BZ\omega_\alpha)$.
\end{lemma}

\begin{proof}
	Let $\gamma:X\to T_0^\vee$ be the restriction map. It induces an isomorphism
	$X/(F-1)X\iso T_0^\vee$.
	Let $I$ be a $\phi$-stable subset of $\Delta$.
	Let $X'_I$ be the set of $\lambda\in X_q$ such that
	\begin{itemize}
		\item given $\alpha{\not\in}I$, we have
			$\langle\lambda,\alpha^\vee\rangle=q-1$ 
		\item given $\alpha\in I$, there is $i\ge 0$ such that
			$\langle\lambda,\phi^i(\alpha^\vee)\rangle\neq q-1$.
	\end{itemize}
	The map $\gamma$ restricts to a bijection $X'_I\iso \gamma(\bigoplus_{\alpha\in I}\BZ\omega_\alpha)$.

	The restriction $\bar{\gamma}_{|X'_I}$ factors as 
	$$X'_I\xrightarrow[\sim]{\gamma_{|X'_I}}
	\gamma(\bigoplus_{\alpha\in I}\BZ\omega_\alpha)\xrightarrow{\gamma_I}Z(G)^\vee$$
	where $\gamma_I$ is given by restricting from $T_0$ to $Z(G)$.

	Let $\zeta\in Z(G)^\vee$. Put $\delta_{I,\zeta}=1$ if $\gamma_I^{-1}(\zeta)\neq\emptyset$ and
	$\delta_{I,\zeta}=0$ otherwise.

	We have 
	$$|\Irr_k(kGb_\zeta)|=\sum_I |\bar{\gamma}^{-1}(\zeta)\cap X'_I|=
	\sum_I |\gamma_I^{-1}(\zeta)|=
	\sum_I \delta_{I,\zeta}|\gamma_I^{-1}(0)|$$
	where $I$ runs over non-empty $\phi$-stable subsets of $\Delta$.
This shows the requested inequality.

	\smallskip
	We have $\delta_{I,\zeta}=1$ if and only if $\zeta\in\bar{\gamma}(\bigoplus_{\alpha\in I}
	\BZ\omega_\alpha)$. Note that $\bar{\gamma}(\omega_\alpha)=\bar{\gamma}(\omega_{\phi(\alpha)})$ for all
	$\alpha$. The equivalence of the lemma follows.
\end{proof}

The tables of [Bki, Lie 4,5,6] show that outside type $A$, there is
a fundamental weight in the root lattice, hence 
$\bigcap_{\alpha\in\Delta^\vee} \bar{\gamma}(\BZ\omega_\alpha)=0$.

Assume $G=\SL_n(q)$ (in which case we put $\eps=1$) or $G=\SU_n(q)$ (in which case we put $\eps=-1$).
We have 
$\bigcap_{\alpha\in\Delta^\vee} \bar{\gamma}(\BZ\omega_\alpha)\neq 0$ if and only if
$n=\ell^r$ for some prime $\ell|(q-\eps)$ and $r\ge 1$.
In that case $Z(G)\simeq\BZ/(\gcd(\ell^r,q-\eps)\BZ)$ and the non-trivial characters of $Z(G)^\vee$
with order $\ell$ are those non-trivial characters $\zeta$ such that $kGb_\zeta$ has the same
number of simple modules as $kGb_0$.

\subsubsection{Equivalences}

We consider the setting above where $G=\SL_{\ell^r}(q)$ or $G=\SU_{\ell^r}(q)$ with $\ell$ a prime
dividing $q-\eps$.

\begin{question}
	Let $\zeta$ be a character of order $\ell$ of $Z(G)^\vee$. What is the relation between
	$kGb_0\mMod$ and $kGb_\zeta\mMod$? Are the blocks
	$kGb_0$ and $kGb_\zeta$ stably equivalent?
\end{question}

Note that if $\BP$ is a proper $F$-stable parabolic subgroup of $\BG$, then the character $\zeta$
of $Z(G)$ extends to $P$. Given $Q$ a non-trivial $p$-subgroup of $G$, there is a proper
$F$-stable parabolic subgroup $\BP$ of $\BG$ such that $N_G(Q)\subset P$ (cf \cite[Remark 6.15]{CaEn2}).
It follows that $\zeta$ extends to $N_G(Q)$. As a consequence (cf Remark \ref{re:blockscenter} below), 
the blocks of proper local subgroups corresponding to $kGb_0$ and $kGb_\zeta$ are isomorphic.
If those local equivalences could be glued (cf \cite[\S 7.3]{Rou1} for a setting for gluing),
then we would obtain a stable equivalence answering positively the question.

\medskip
	When $\ell^r=2$ and $\eps=1$,
	the question has a positive answer. The blocks are actually derived
	equivalent, since they are both derived equivalent to the corresponding blocks of $kB_0$
	(cf \S\ref{se:abelianp}), and those blocks are isomorphic. Cf also \cite{Wo1} for a direct
	construction as a composition of perverse equivalences.

	When $\ell^r=3$ and $\eps=-1$, the question has also a positive answer since the blocks
	are both stably equivalent to the corresponding blocks of $kB_0$
	(cf \S\ref{se:rank1}), and those blocks are isomorphic.

\begin{rem}
Note that we do not expect the blocks to be derived equivalent in general. Consider the case
$G=\SL_3(4)$. Let us show that there are no perfect isometries between the principal $2$-block and
a non-principal $2$-block of $\mathrm{SL}_3(4)$. In particular, those blocks are not derived equivalent
over $\mathcal{O}$.

\smallskip
The principal block of $G$ has $9$ irreducible characters: $1$, $20$, $35$, $\overline{35}$, $35'$, $45$, $\overline{45}$, $63$ and $\overline{63}$. Let $b$ be one of the non-principal block idempotents with positive defect. The irreducible characters of $KGb$ are those with central character a given primitive cubic root of unity. They are $15_1$, $15_2$, $15_3$, $21$, $45$, $\overline{45}$, $63$, $\overline{63}$ and $84$. 

Let $I$ be a perfect isometry from $\CO Gb_0$ to $\CO Gb$. 
Let $\eta=\sum_{\chi\in\mathrm{Irr}(KGb_0)}\chi\otimes I(\chi)$. Let $g$ be an element of order $5$ and $h$ an element of even order. We have 
$$0=\eta(g,h)=I(1)(h)+\alpha I(63)(h)+\bar{\alpha}I(\overline{63})(h)$$
with $\alpha=\frac{1-\sqrt{5}}{2}$ and $\bar{\alpha}=\frac{1+\sqrt{5}}{2}$. It follows that $I(63)(h)=I(\overline{63})(h)=-I(1)(h)$.  Taking $h$ an involution, we deduce that $\{I(1),I(63), I(\overline{63}\}\subset\{\pm 15_1,\pm 15_2,\pm 15_3,\pm 63,\pm\overline{63}\}$. Taking $h$ of order $4$, we obtain that $\{I(1),I(63), I(\overline{63}\}\subset\{\pm 21,\pm 45,\pm \overline{45},\pm 63,\pm\overline{63}\}$.  This is a contradiction.
\end{rem}

\begin{rem}
	\label{re:blockscenter}
	Let $G$ be a finite group and $\zeta$ a character of $Z(G)_{p'}$, the largest $p'$-subgroup of
	$Z(G)$. Let
	$$e_\zeta=\frac{1}{|Z(G)_{p'}|}\sum_{z\in Z(G)_{p'}}\zeta(z)^{-1}z$$
	be the associated 
idempotent.

\smallskip
If $\zeta$ extends to a character $\tilde{\zeta}$ of $G$, then there is an isomorphism of algebras 
$$e_1 kG\iso e_\zeta kG,\ \sum_{g\in G}a_g g\mapsto \sum_{g\in G}\tilde{\zeta}(g)^{-1}a_gg.$$
In general, the algebras $e_1kG$ and $e_\zeta kG$ can have rather different module categories and
invariants.
When the extension assumption holds for all $p$-local subgroups, we could hope
that the algebras are at least stably equivalent. More precisely, 
	assume that given any non-trivial $p$-subgroup $Q$ of $G$, the character $\zeta$ extends
	to $N_G(Q)$. Are the algebras $e_1kG$ and $e_\zeta kG$ stably equivalent?

\smallskip
	One can ask similar questions for two linear characters of $Z(G)$. For example, let $d$ be
	an integer prime to the order of $\zeta$. What is the relation between $e_\zeta kG$ and
	$e_{\zeta^d}kG$? If $\zeta^d$ can be obtained from $\zeta$ by applying a field 
	automorphism of $k$, then the rings $e_\zeta kG$ and
	$e_{\zeta^d}kG$ are isomorphic, but they need not be isomorphic as $k$-algebras  (nor even derived
	equivalent), as shown by Benson and Kessar \cite[Example 5.1]{BeKe}. In their examples,
	$O_p(G)\neq 1$.
\end{rem}

\section{Non-defining characteristic}
\label{se:nondef}

We consider a connected
reductive algebraic group $\BG$ with a Steinberg endomorphism $F$ as in \S \ref{se:Steinberg}.
We fix a prime $\ell$ distinct from $p$. We will be discussing mod-$\ell$ representations of $G$.
 We fix $K$ a finite extension of $\BQ_{\ell}$ containing all
$|G|$-th roots of unity and denote by $\CO$ its ring of integers and by $k$ its residue field.

In \S\ref{se:DL} and \S\ref{se:char0}, we recall constructions and results of Deligne-Lusztig and Lusztig
\cite{DeLu1,DeLu2,Lu1,Lu2}.
\subsection{Deligne-Lusztig varieties}
\label{se:DL}

\subsubsection{Definition}
Consider the Lang covering $\CL:\BG\to\BG,\ g\mapsto g^{-1}F(g)$. This is a surjective \'etale Galois
morphism, with Galois group $G$.

\smallskip
Let $\BP$ be a parabolic subgroup of $\BG$ and let $\BU$ be its unipotent radical. Assume there is an
$F$-stable
Levi subgroup $\BL$ with $\BP=\BU\rtimes\BL$.
The associated Deligne-Lusztig variety is $\CL^{-1}(F(\BU))$.
It has a free left (resp. right) action of
$G$ (resp. $L$) by multiplication. We can also consider its quotient $Y_{\BU}=\CL^{-1}(F(\BU))/(\BU\cap F(\BU))$,
which has the same $\ell$-adic cohomology. One can consider further the variety
$X_{\BU}=Y_{\BU}/L$. The varieties $Y_{\BU}$ and $X_{\BU}$ are smooth.

\begin{rem}
	When $\BP$ is $F$-stable, then $Y_\BU=G/U$ is a finite set.
\end{rem}

\subsubsection{Case of tori}
\label{se:casetori}
A particular role is played by Deligne-Lusztig varieties associated to tori. Let us give another model
for those.
Fix $\BB_0$ an $F$-stable Borel subgroup of $\BG$ and $\BT_0$ an $F$-stable maximal torus contained
in $\BB_0$. Let $\BU_0$ be the unipotent radical of $\BB_0$ and let $W=N_\BG(\BT_0)/\BT_0$. 

\smallskip
Let $\CB$ be the variety of Borel subgroups of $\BG$. There is a decomposition
$\CB\times\CB=\coprod_{w\in W}\CO(w)$ into orbits for the diagonal action of $\BG$, where
$\CO(w)$ is the orbit containing $(\BB_0,w\BB_0 w^{-1})$.

Given $w\in W$, we put 
$$X(w)=\{\BB\in\CB\ |\ (\BB,F(\BB))\in\CO(w)\}.$$

\smallskip
Let $\dot{w}\in N_{\BG}(\BT_0)$ with image $w\in W$. We put
$$Y(\dot{w})=\{g\BU_0\in\BG/\BU_0\ |\ g^{-1}F(g)\in \BU_0\dot{w}\BU_0\}.$$
There is a left action of $G$ on $X(w)$ and $Y(\dot{w})$ by left multiplication and a right action of
$\BT_0^{wF}$ on $Y(\dot{w})$ by right multiplication. The map $g\BU_0\mapsto g\BB_0g^{-1}$ induces
an isomorphism of $G$-varieties $G\setminus Y(\dot{w})\iso X(w)$. The varieties $X(w)$ and $Y(\dot{w})$
have pure dimension $l(w)$, the length of $w$ (cf \S\ref{se:real}).

\smallskip
Let $h\in\BG$ such that $h^{-1}F(h)=\dot{w}$.
The maximal torus $\BT=h\BT_0 h^{-1}$ is $F$-stable and 
the isomorphism $\BT_0\iso\BT,\ t\mapsto hth^{-1}$ restricts to an isomorphism
$\BT_0^{wF}\iso \BT^F$. 

There is a commutative diagram 
$$\xymatrix{
	Y_{h\BU_0 h^{-1}}\ar[rrr]^{g\mapsto gh\BU_0}_\sim 
	\ar[d]_{g(\BU\cap F(\BU))\mapsto g(\BU\cap F(\BU))L} &&&
Y(\dot{w}) \ar[d]^{g\BU_0\mapsto g\BB_0g^{-1}}\\
X_{h\BU_0 h^{-1}}\ar[rrr]_{g\mapsto gh\BB_0h^{-1}g^{-1}}^\sim  &&& X(w)}$$
where the horizontal maps are $G$-equivariant isomorphisms and the top horizontal map is equivariant
for the right action of $T$, via its identification with $\BT_0^{wF}$ above.

\medskip
Two elements $w$ and $w'$ of $W$ are {\em $F$-conjugate} if there is $v\in W$ such that
$w'=v^{-1}wF(v)$.

\smallskip
The construction of $\BT$ from $w$ induces a bijection from the
set of $F$-conjugacy classes of $W$ to the set of $G$-conjugacy classes of $F$-stable maximal
tori of $\BG$.

\begin{rem}
	The varieties $X(w)$ are known to be affine in many cases, but it is not known if they are affine
	in general.
The affinity is known when $q$ is larger than the Coxeter number of $\BG$ 
	\cite[Theorem 9.7]{DeLu1} and when $w$ has minimal length in its $F$-conjugacy class
	(\cite[Corollary 2.8]{Lu1}, \cite[\S 5]{OrRa}, \cite[Theorem 1.3]{He} and \cite{BoRou3}).
\end{rem}

\subsubsection{Endomorphisms}
\label{se:endomorphisms}
We follow \cite{BrouMi2}, inspired by an earlier construction of Lusztig \cite[pp. 24--25]{Lu2}.

Given $w_1,\ldots,w_r\in W$, let 
$$X_F(w_1,\ldots,w_r)=\{(\BB_1,\ldots,\BB_r)\in\CB^r\ |\
(\BB_i,\BB_{i+1})\in\CO(w_i),\ 1\le i<r \text{ and }(\BB_r,F(\BB_1))\in\CO(w_r)\}.$$
The variety $X_F(w_1,\ldots,w_r)$ depends only on the element $b=\lambda(w_1)\cdots\lambda(w_r)$
of $B_W^+$ (cf \S\ref{se:real} for the notations),
up to canonical isomorphism \cite{De} and we denote it by $X_F(b)$.

There is an action of $\phi$ on $B_W$ by $b_s\mapsto b_{\phi(s)}$ for $s\in S$.
Note that given $n>0$, we have a morphism 
\begin{align*}
	\iota_n:X_F(b)&\to X_{F^n}(b\phi(b)\cdots \phi^{n-1}(b))\\
	(\BB_1,\ldots,\BB_r)&\mapsto
(\BB_1,\ldots,\BB_r,F(\BB_1),\ldots,F(\BB_r),\ldots,F^{n-1}(\BB_1),\ldots,F^{n-1}(\BB_r)).
\end{align*}

Given $0\le i\le r$, the morphism 
\begin{align*}
	D_{w_1,\ldots,w_i}:X_F(w_1,\ldots,w_r)&\to X_F(w_{i+1},\ldots,w_r,\phi(w_1),\ldots,\phi(w_i))\\
	(\BB_1,\ldots,\BB_r)&\mapsto (\BB_{i+1},\ldots,\BB_r,F(\BB_1),\ldots,F(\BB_i))
\end{align*}
is purely inseparable.
Given $b',b''\in B_W^+$, this provides a morphism $D_{b'}:X_F(b'b'')\to X_F(b''\phi(b'))$.

\medskip
Let $b\in B_W^+$ such that $(b\phi)^d=\pi^r\phi^d$ for some $d,r>0$. We define an action of
$C_{B_W^+}(b\phi)$ on $X_F(b)$ as follows.

Let $b'\in C_{B_W^+}(b\phi)$. There are $t>0$ and $b''\in B_W^+$ such that $\pi^{rt}=bb'b''$.
We identify $X_F(b)$ with a subvariety of $X_{F^{dt}}(\pi^{rt})$ using the embedding $\iota_{dt}$.
The endomorphism $D_{b'}$ of $X_{F^{dt}}(\pi^{rt})$ preserves $X_F(b)$

\medskip
The constructions above extend to the varieties $Y$.
The composite morphism $B_W\xrightarrow{\can} W=N_{\BG}(\BT_0)/\BT_0$ lifts to a morphism
$\sigma:B_W\to N_{\BG}(\BT_0)$. Given $w_1,\ldots,w_r\in W$, we have a variety
$$Y_F(w_1,\ldots,w_r)=\{(g_1\BU_0,\ldots,g_r\BU_0)\in (\BG/\BU_0)^r\ |\
\begin{cases}
	g_i^{-1}g_{i+1}\in\BU_0\sigma(\lambda(w_i))\BU_0  \text{ for } 1\le i<r \\
g_r^{-1}F(g_1)\in \BU_0\sigma(\lambda(w_r))\BU_0
\end{cases}
\}.$$
It has a left action of $G$ by diagonal left multiplication and a right action of $\BT_0^{w_1\cdots w_r F}$
by diagonal right multiplication.
We have a $G$-equivariant morphism corresponding to the quotient by $\BT_0^{w_1\cdots w_r F}$
$$Y_F(w_1,\ldots,w_r)\to X_F(w_1,\ldots,w_r),\
(g_1\BU_0,\ldots,g_r\BU_0)\mapsto (g_1\BB_0g_1^{-1},\ldots,g_r\BB_0g_r^{-1}).$$

The variety $Y_F(w_1,\ldots,w_r)$ depends only on the element $b=\lambda(w_1)\cdots\lambda(w_r)$
of $B_W^+$, up to canonical isomorphism and we denote it by $Y_F(b)$. 
Given $b',b''\in B_W^+$, we obtain a morphism $D_{b'}:Y_F(b'b'')\to Y_F(b''\phi(b'))$ and we have
an action of $C_{B_W^+}(b\phi)$ on $Y_F(b)$ compatible with the action on $X_F(b)$.

\subsubsection{Deligne-Lusztig functors}

Let $\BP$ be a parabolic subgroup of $\BG$ with unipotent radical $\BU$ and an $F$-stable Levi
complement $\BL$. The complex $\Lambda_c(Y_\BU,\BZ_{\ell})$ of $\ell$-adic cohomology with compact
support of $Y_{\BU}$ 
\cite{Ric1,Rou2} is a bounded complex of $\ell$-permutation $\BZ_{\ell}(G\times L^\opp)$-modules. It is 
well defined up to homotopy. Its cohomology groups are the $H^i_c(Y_\BU,\BZ_{\ell})$.

\smallskip
Given $R$ a commutative $\BZ_{\ell}$-algebra, we put $\Lambda_c(Y_\BU,R)=\Lambda_c(Y_\BU,\BZ_\ell)
\otimes_{\BZ_\ell}R$. We obtain a functor
$$\CR_{\BL\subset\BP}^{\BG}=\Lambda_c(Y_\BU,R)\otimes_{RL}-:D^b(RL\mMod)\to D^b(RG\mMod).$$

\smallskip
When $R=\bar{\BQ}_\ell$, the functor $\CR_{\BL\subset\BP}^{\BG}$ induces a morphism
$$R_{\BL\subset\BP}^{\BG}:G_0(\bar{\BQ}_\ell L)\to G_0(\bar{\BQ}_\ell G).$$

Note that this morphism is expected to depend only on $\BL$, and not on $\BP$. This is known to hold
except possibly when $q=2$, the parabolic subgroup $\BP$ is not $F$-stable
and the Dynkin diagram of $\BG$ contains a subdiagram of type $E_6$
(cf \cite[Corollary 4.3]{DeLu1} and \cite{DeLu2,BoMi}).

\begin{rem}
	When $\BP$ is $F$-stable, then $\CR_{\BL\subset\BP}=R[G/U]\otimes_{RL}-$ is the Harish-Chandra
	induction functor.
\end{rem}

We denote by ${^*R}_{\BL\subset\BP}^\BG: G_0(\bar{\BQ}_\ell G)\to G_0(\bar{\BQ}_\ell L)$ the adjoint
of $R_{\BL\subset\BP}^{\BG}$.

\subsection{Characteristic $0$ representations}
\label{se:char0}
\subsubsection{Tori and characters}
\label{se:torichar}
Let $\BT$ be an $F$-stable maximal torus of $\BG$.
Fix $M$ a positive integer multiple of $\delta$ such that $(wF)^M(t)=t^{q^M}$ for all $t\in 
\BT$ and $w\in W$ (cf \S\ref{se:Steinberg} for the definitions of $\delta$ and $q$).
Let $\zeta$ (resp. $\xi$) be a root of unity of order $q^M-1$ of $\bar{\BF}_p$
(resp. $\bar{\BQ}_\ell$).

The morphism 
$$N:Y(\BT)\to T,\ y\mapsto y(\zeta) F(y(\zeta))\cdots
F^{M-1}(y(\zeta))$$
is surjective and induces an isomorphism $Y(\BT)/\bigl((F-1)(Y(\BT)\bigr)\iso T$.

The morphism
$$X(\BT)\to \Hom(Y(\BT),\bar{\BQ}_{\ell}^\times),\  \chi\mapsto (y\mapsto
\xi^{\langle \chi, y+F(y)+\cdots+F^{M-1}(y)\rangle})$$
factors through
$\Hom(N,\bar{\BQ}_{\ell}^\times)$ and gives a surjective morphism
$X(\BT)\to \Hom(T,\bar{\BQ}_{\ell}^\times)$. That induces an isomorphism
$X(\BT)/\bigl((F-1)(X(\BT)\bigr)\iso \Irr_{\bar{\BQ}_\ell}(T)$.

\subsubsection{Tori and dual groups}
\label{se:toridual}
Let $(\BG^*,\BT_0^*,F^*)$ be a triple dual to $(\BG,\BT_0,F)$: the group $\BG^*$ is a Langlands dual of
$\BG$, there is a given isomorphism $X(\BT_0^*)\iso Y(\BT_0)$, and $F^*$ is a Steinberg endomorphism of
$\BG^*$ stabilizing $\BT_0^*$ and dual to $F$. Furthemore, there
is a given isomorphism $W^*=N_{\BG^*}(\BT^*_0)/\BT^*_0\iso W=N_{\BG}(\BT_0)/\BT_0$ and we identify those
groups. Note that the action
of $F^*$ on $W^*$ corresponds to the action of $F^{-1}$ on $W$.

\smallskip
Let $\BT$ be an $F$-stable maximal torus of $\BG$. It corresponds to an $F$-conjugacy class
$(w)$ of $W$ (cf \S \ref{se:casetori}). We denote by $\BT^*$ an $F^*$-stable maximal torus of $\BG^*$
whose $(\BG^*)^{F^*}$-conjugacy class is given by the $F^*$-conjugacy class $(w^{-1})$.
Furthermore, the identification of $\BT_0^*$ with the dual of $\BT_0$ provides an isomorphism
between $\BT^*$ and the dual of $\BT$, and that isomorphism is well-defined up to the action
of $(N_{\BG}(\BT)/\BT)^F$. Via the constructions of \S\ref{se:torichar}, this gives an isomorphism
$\Irr_{\bar{\BQ}_\ell}(T)\iso (\BT^*)^{F^*}$.

This construction provides a bijection from the set of $G$-conjugacy classes of pairs $(\BT,\theta)$,
where $\BT$ is an $F$-stable maximal torus of $\BG$ and $\theta\in\Irr_{\bar{\BQ}_\ell}(\BT^F)$
to the set of $(\BG^*)^{F^*}$-conjugacy classes of pairs $(\BT^*,s)$ where $\BT^*$ is an 
$F^*$-stable maximal torus of $\BG^*$ and $s\in (\BT^*)^{F^*}$.

\subsubsection{Jordan-Lusztig decomposition}
Let us recall the Jordan decomposition of conjugacy classes. An element $g\in \BG$ can be decomposed
uniquely as $g=tu$ where $t$ is semi-simple, $u$ is unipotent and $ut=tu$. 
Denote by $\mathrm{Cl}(G)$ (resp. $\mathrm{Cl}_{ss}(G)$, $\mathrm{Cl}_{unip}(G)$) the set of conjugacy classes of
elements (resp. semi-simple, unipotent elements) of $G$.

The Jordan decomposition induces a bijection
$$\mathrm{Cl}(G)\iso \coprod_{(t)\in \mathrm{Cl}_{ss}(G)}\mathrm{Cl}_{unip}(C_G(t))$$
where
$t$ runs over conjugacy classes of semi-simple elements of $G$.

\smallskip
Given $(s)$ a conjugacy class of semi-simple elements of $(\BG^*)^{F^*}$, we denote by
$\Irr_{\bar{\BQ}_\ell}(G,(s))$ the set of irreducible representations of $G$ that occur in the $\theta$-isotypic
component of $H^*_c(Y_\BU,\bar{\BQ}_\ell)$ for some Borel subgroup of $\BG$ with unipotent radical $\BU$ and
containing an $F$-stable maximal torus $\BT$ and $\theta\in\Irr_{\bar{\BQ}_\ell}(T)$ such that
$(\BT,\theta)$ corresponds to
$(\BT^*,s)$ by the bijection of \S\ref{se:toridual} for some
$F^*$-stable maximal torus $\BT^*$ of $\BG^*$ containing $s$.

The {\em unipotent representations} of $G$ are those in $\Irr_{\bar{\BQ}_\ell}(G,1)$. They are the irreducible 
representations of $G$ that occur in $H_c^*(X(w),\bar{\BQ}_\ell)$ for some $w\in W$.

\smallskip

We have the Deligne-Lusztig decomposition (cf \cite[Theorem 8.24]{CaEn2})
$$\Irr_{\bar{\BQ}_\ell}(G)=\coprod_{(s)\in\mathrm{Cl}_{ss}(\BG^*)^{F^*}}\Irr_{\bar{\BQ}_\ell}(G,(s))$$
where $(s)$ runs over conjugacy classes of semi-simple elements of $(\BG^*)^{F^*}$.

\smallskip

Let $s$ be a semi-simple element of $(\BG^*)^{F^*}$. When $C_{\BG^*}(s)$ is connected, let
$(C_{\BG^*}(s)^*,F)$ be dual to $(C_{\BG^*}(s),F^*)$. Note that 
$C_{\BG^*}(s)^*$ need not occur as a subgroup of $\BG$.

Lusztig constructed a bijection (cf \cite[Theorem 15.8]{CaEn2})
$$\Irr_{\bar{\BQ}_\ell}(\bigl(C_{\BG^*}(s)^*\bigr)^F,1)\iso \Irr_{\bar{\BQ}_\ell}(G,(s)).$$

When $C_{\BG^*}(s)$ is a Levi subgroup of $\BG^*$, then $C_{\BG^*}(s)^*$ can be realized as an
$F$-stable Levi subgroup $\BL$ of $\BG$ and the bijection is given by
\begin{equation}
	\label{eq:JordanDL}
	\rho\mapsto \pm R_{\BL\subset\BP}^{\BG}(\rho\otimes\eta)
\end{equation}
where $\eta$ is the one-dimensional representation of $L$ corresponding by duality to 
$s\in Z(\BL^*)^{F^*}$ and $\BP$ is a parabolic subgroup of $\BG$ with Levi complement $\BL$.

\smallskip
When $Z(\BG)$ is connected, one obtains 
the Jordan-Lusztig decomposition of characters
$$\Irr_{\bar{\BQ}_\ell}(G)\iso \coprod_{(s)\in\mathrm{Cl}_{ss}(\BG^*)^{F^*}}\Irr_{\bar{\BQ}_\ell}(\bigl(C_{\BG^*}(s)^*\bigr)^F,1).$$
If $Z(\BG)$ is connected, then $\mathrm{Cl}_{ss}(\BG^*)^{F^*}$ can be replaced by
$\mathrm{Cl}_{ss}((\BG^*)^{F^*})$.

\subsubsection{Unipotent representations}
\label{se:unipotentrep}

Lusztig constructed a parametrization of simple
unipotent representations of $G$ by a combinatorially defined
set $\CU(W,\phi)$
that depends only on the Weyl group $W$ and on the finite order automorphism $\phi$ induced by $F$ of the
reflection representation of $W$. The degrees of the irreducible unipotent representations are
polynomials in $q$ (the {\em generic degrees}).
Lusztig also defined a partition of $\CU(W,\phi)$ into {\em families}, and a partial order on
the set of families.

\medskip
When $\BG=\GL_n$, the simple unipotent representations are parametrized by partitions of $n$. When
$G=\GL_n(q)$, the simple unipotent representations are the components
of $\Ind_{B_0}^G\bar{\BQ}_\ell$.

\subsection{Modular representations}
\label{se:modular}

\subsubsection{Blocks and Lusztig series}
Let $t$ be a semi-simple element of $(\BG^*)^{F^*}$ of order prime to $\ell$. 
We put 
$$e_{(t)}=e_{(t)}^G=
\sum_{\substack{(s)\in \mathrm{Cl}_{\ell}(C_{\BG^*}(t)^{F^*})\\ \chi\in\Irr_{\bar{\BQ}_\ell}(G,
(st))}}e_\chi$$
where $\mathrm{Cl}_\ell$ denotes the set of conjugacy classes of $\ell$-elements.

This idempotent of $Z(KG)$ is actually in $Z(\CO G)$ \cite{BrouMi1}, hence it is a sum of (orthogonal) block
idempotents. In other terms, 
$\bigcup_{(s)\in \mathrm{Cl}_{\ell}(C_{\BG^*}(t)^{F^*})}\Irr_{\bar{\BQ}_\ell}(G,(st))$ is a union
of characters in blocks.

A {\em unipotent block} is a block $kGb$ such that $be_{(1)}=b$.

\smallskip
Given $\BB$ a Borel subgroup of $\BG$ containing an $F$-stable maximal torus $\BT$, with unipotent radical
$\BU$ and given $\theta\in\Irr(T)_{\ell'}$ such that the pair $(\BT,\theta)$ corresponds to a pair
$(\BT^*,t)$, then $\Lambda_c(Y_\BU,\CO)e_\theta$ is an object of $\CO Ge_{(t)}\mperf$. Furthermore, those
complexes (for varying $\BB$, $\BT$ and $\theta$) generate $\CO Ge_{(t)}\mperf$ (the smallest full
thick triangulated subcategory containing those is the whole category) \cite[Theorem A']{BoRou1}.

There
is a similar statement for derived categories when all elementary abelian $\ell$-subgroups of $G$ are
contained in tori \cite[Theorem 1.2]{BoDaRou}.

\subsubsection{Jordan decomposition}
Brou\'e conjectured \cite{Brou2} a modular version of (\ref{eq:JordanDL}): assume
$C_{\BG^*}(t)$ is a Levi subgroup of $\BG^*$, with corresponding dual an $F$-stable Levi
subgroup $\BL$ in $\BG$, and let $\eta$ be dual to $t$.
Let $\BP$ be a parabolic subgroup of $\BG$ with unipotent radical $\BU$ and Levi complement $\BL$.
Then
$H^{\dim Y_\BU}(Y_\BU,\CO)\otimes\eta$ induces a Morita equivalence between $\CO Ge_{(t)}^G$ and
$\CO Le_{(1)}^L$. This was proven in \cite{Brou2} when $\BL$ is a torus, while
the general case is \cite[Theorem 11.8]{BoRou1}. There is an extension of that result to the case where
$C_{\BG^*}(t)^\circ$ is a Levi subgroup \cite{BoDaRou}.

\begin{rem}
The geometric approach has not enabled us to relate isolated blocks,
i.e. corresponding to a semi-simple $\ell'$-element 
	$t$ of $(\BG^*)^{F^*}$ such that $C_{\BG^*}(t)^\circ$
is not contained in a proper Levi subgroup, to unipotent blocks. It is conjectured though that any block
is Morita equivalent to a unipotent block, for a possibly non-connected group.
\end{rem}

\subsubsection{Unipotent blocks}

Assume $p$ and $\ell$ are good for $\BG$ and $\ell{\not|}|Z(\BG)/Z(\BG)^\circ|\cdot
|Z(\BG^*)/Z(\BG^*)^\circ|$. We also assume for the remainder of \S\ref{se:modular} that $F$ is a Frobenius
endomorphism.

\smallskip
There is a "$d$-Harish Chandra" parametrization of blocks of $kG$ \cite{CaEn1,BrouMalMi}.

An $F$-stable Levi subgroup of $\BG$ is {\em $d$-split} if it is the centralizer of a $\Phi_d$-subgroup
of $\BG$. A simple unipotent representation $\rho$ of $G$ is said to be {\em $d$-cuspidal} if
${^*R_{\BL\subset\BP}^\BG}(\rho)=0$ for all proper $d$-split Levi subgroups $\BL$ of $\BG$ with
$\BP$ a parabolic subgroup of $\BG$ with Levi complement $\BL$.

\smallskip
There is a parametrization of the set of unipotent blocks by 
the set $G$-conjugacy classes of pairs $(\BL,\lambda)$ where $\BL$ is a
$d$-split Levi subgroup of $\BG$ and $\lambda$ is a $d$-cuspidal unipotent character of $G$:
the simple unipotent representations in the block corresponding to $(\BL,\lambda)$ are those
that occur in $R_{\BL\subset\BP}^\BG(\lambda)$ for some $\BP$.

The parametrization of classes of pairs $(\BL,\lambda)$ above depends only on $(W,\phi)$ and $d$, and
the corresponding subset of $\CU(W,\phi)$ depends only on $d$ \cite{BrouMalMi}.

\subsubsection{Unipotent decomposition matrices}
\label{se:decmat}
Fix a total order on $\Irr_{\bar{\BQ}_\ell}(G,1)$ such that if $\rho_i$ is in the family $\CF_i$ for
$i\in\{1,2\}$ and $\CF_1<\CF_2$, then $\rho_1<\rho_2$ (cf \S\ref{se:unipotentrep}).

\smallskip
The following result was conjectured by Geck \cite{Ge1} and proven by \cite{BruDuTa}, following earlier work
on basic sets \cite{Ge2,GeHi} and proofs for $\GL_n(q)$ in \cite{DipJa}, for $\mathrm{GU}_n(q)$ in
\cite{Ge2} and for classical groups
and certain $\ell$ (linear primes, for which the blocks are related to blocks of $\GL_n(q)$) 
in \cite{GruHi}. 

\begin{thm}
	\label{th:triangular}
	There is a (unique) bijection $\beta:\Irr_{\bar{\BQ}_\ell}(G,1)\iso \Irr(kGe_1)$ such that
	$\mathrm{dec}([\rho])\in [\beta(\rho)]+\sum_{\rho'>\rho}\BZ_{\ge 0}[\beta(\rho')]$ for any
	$\rho\in\Irr_{\bar{\BQ}_\ell}(G,1)$.
\end{thm}

The theorem above together with Lusztig's work (\S \ref{se:unipotentrep}) provides a 
parametrization of the set of simple
$kGe_1$-modules by a set that depends only on $(W,\phi)$.

\smallskip
It is conjectured that, given $W$ and $d$ (the order of $q$ in $\BF_\ell^\times$),
for $\ell$ large enough, the square part of the
decomposition matrix involving unipotent representations depends only on $(W,\phi)$ and $d$,
i.e., it is independent of $\ell$ and $q$. This
is known for $\GL_n(q)$ \cite{DipJa} and for linear primes and classical groups \cite{GruHi}.

The determination of this "generic" square matrix is a major open problem in the study of decomposition
matrices for finite groups of Lie type in non-defining characteristic. The recent \cite{DuMal2} 
provides a number of new decomposition matrices for groups of low rank.

\medskip
Assume $\BG$ is split.
The algebra $\End_{\CO G}(\Ind_{B_0}^G\CO)$ is isomorphic to the Hecke algebra of $W$ over $\CO$
(cf \S\ref{se:real}),
specialized at $x=q$. The decomposition matrix of that specialized Hecke algebra is equal to
the submatrix of the decomposition matrix with rows parametrized by simple modules that are
direct summands of $\Ind_{B_0}^G K$ (principal series representations)
and columns by simple modules that are quotients of $\Ind_{B_0}^Gk$ \cite{Dip}.
The former depends only on $d$, if $\ell$ is large enough, as it is the same as the one for the Hecke
algebra at $x$ a primitive $d$-th root of unity, over $\BC$ \cite{Ge3}.
This shows the genericity property for a small submatrix.
Similar considerations can be used to prove genericity properties for small submatrices along the diagonal
corresponding to various Harish-Chandra series using relative Hecke algebras.

\medskip
Let $a=\nu_\ell(\Phi_d(q))$.
Theorem \ref{th:triangular} asserts that the decomposition matrix has the following shape. Here,
the gray entries are on rows corresponding to principal series irreducible characters and columns
corresponding to modular simple representations that are quotients of $\Ind_{B_0}^Gk$. The gray entries
give the decomposition matrix of the Hecke algebra.

$$
\begin{tikzpicture}[
		squarednode/.style={rectangle, fill=gray!25, very thick, minimum size=5mm},]
		
	\draw (0,0) -- (0,-1) -- (6,-1) -- (6,0);
	\draw[dashed] (0,0) -- (0,1);
	\draw[dashed] (6,0) -- (6,1);
	\draw (0,1) -- (0,8) -- (6,8) -- (6,1);
	\draw[line width=1mm] (0.1,2) -- (0.1,7.9) -- (5.9,7.9) -- (5.9,2) -- (0.1,2);
%	\draw[line width=1mm,color=gray] (0.2,4.1) -- (0.2, 7.8) -- (1.8, 7.8) -- (1.8,4.1) -- (0.2,4.1);
	\draw (0.7,7.2) node[squarednode, anchor=south east]{$1$};
	\draw[dashed] (1.1,6.8) node[anchor=south east]{$1$} -- (1.3,6.6) node[anchor=north west]{$1$};
	\draw (2.2,5.6) node[squarednode, anchor=south east]{$1$};
	\draw (0.7,5.6) node[squarednode, anchor=south east]{};
	\draw (0.7,4.6) node[squarednode, anchor=south east]{};
	\draw (2.2,4.6) node[squarednode, anchor=south east]{};
	\draw[dashed] (2.6,5.2) node[anchor=south east]{$1$} -- (2.9,4.9) node[anchor=north west]{};
	\draw[dashed] (3.0,4.8) node[anchor=south east]{} -- (3.2,4.6) node[anchor=north west]{$1$};
	\draw (4.1,3.6) node[squarednode, anchor=south east]{$1$};
	\draw (0.7,3.6) node[squarednode, anchor=south east]{};
	\draw (2.2,3.6) node[squarednode, anchor=south east]{};
	\draw (0.7,2.5) node[squarednode, anchor=south east]{};
	\draw (2.2,2.5) node[squarednode, anchor=south east]{};
	\draw (4.1,2.5) node[squarednode, anchor=south east]{};
	\draw[dashed] (4.5,3.1) node[anchor=south east]{$1$} -- (5.4,2.6) node[anchor=north west]{$1$};
	\draw (2.8,3.2) node{$*$};
	\draw (2.8,0.5) node{$*$};
%	\draw[decorate, decoration = {calligraphic brace, amplitude=10pt}] (-0.2,4.1) -- (-0.2,7.8);
	\draw (-2.2,6) node[squarednode]{principal series};
	\draw[decorate, decoration = {calligraphic brace, amplitude=10pt}] (-0.2,-0.9) -- (-0.2,1.9);
	\draw (-2.7,0.5) node{\# rows depends on $\ell^a$};
%	\draw[-stealth] (1.0,8.4) node[squarednode, anchor=south]{Hecke} -- (1.0,7.4);
	\draw (1.0,8.4) node[squarednode, anchor=south]{Hecke};
	\draw[decorate, decoration = {calligraphic brace, amplitude=10pt}] (6.2,7.9) -- (6.2,2.0);
	\draw (7.5,5) node{unipotent};
\end{tikzpicture}$$

\medskip
Assume $G=\GL_n(q)$. Let $A=\CO G/\bigl(\CO G\cap (\bigoplus_{\chi{\not\in}\Irr_{\bar{\BQ}_\ell}(G,1)}
e_\chi KG)\bigr)$. The simple unipotent representations of $kG$ are the same as the simple $A$-modules.
The algebra $A$ is Morita equivalent to the $q$-Schur algebra of $\GS_n$ over $\CO$ specialized
at $x=q$ \cite{DipJa,Ta}.

The $q$-Schur algebra is the endomorphism ring of the direct sum of induced trivial modules
from Hecke algebras of all standard parabolic subgroups.
For $\ell$ large enough, its decomposition matrix is the same as the one obtained 
for $x$ a primitive $d$-the root of unity over $\BC$.
So, the square part of the unipotent decomposition matrix of $\GL_n(q)$ coincides with the decomposition
matrix of the $q$-Schur algebra in characteristic $0$, at a primitive $d$-th root of unity. One deduces
the genericity property for decomposition matrices of $\GL_n(q)$. Furthermore, this matrix has
a description in terms of the combinatorics of the canonical basis
of the Fock space for the quantum group of $\Gsl_d$ \cite{LaLeTh,LeTh,Ar}.

\begin{rem}
One can define analogs of $q$-Schur algebras by generalizing the construction to other types of groups,
but they do not seem to have good descriptions nor good properties like quasi-heredity, except under
particular assumptions making the category of representations look like the one for general linear groups
(for example, classical groups and linear primes). The case of unipotent blocks with cyclic defect, fully
understood now \cite{CrDuRou}, shows already the substantial complications related to the presence of
cuspidal representations. We propose in \S \ref{se:deg} to take a limit $q\to 1$ in the $\ell$-adic
topology, which makes $q\to \infty$ in the real topology.
\end{rem}

\subsection{Brou\'e's conjecture}
\subsubsection{General version}
\label{se:BroueconjP}

Let $b$ be a block idempotent of $\CO G$, $D$ a defect group. Assume $D$ is abelian and
$\BL=C_{\BG}(D)$ is a Levi subgroup of $\BG$. Let $b_D\in \CO L$ be the Brauer correspondent of $b$
and let $b'_D$ be a block idempotent of $\CO L$ with $b'_Db_D=b'_D$.

Given $\BP$ a parabolic subgroup of $\BG$ with unipotent radical $\BU$ and
Levi complement $\BL$, there is a complex of $(\CO G,\CO L)$-bimodules
	$\Lambda_c(Y_{\BU},\CO)b'_D$.

\begin{conj}
	\label{co:Broue2}
	There is a choice of $\BP$ and an extension
	of the right action of $C_G(D)$ on $\Lambda_c(Y_{\BU},\CO)b'_D$ to an action of
	$N_G(D,b'_D)$ such that $\Lambda_c(Y_{\BU},\CO)b'_D$ induces a Rickard equivalence between
	$b\CO G$ and $b'_D\CO N_G(D,b'_D)$.
\end{conj}

We refer to \S\ref{se:equivalences} for the notion of Rickard equivalences.

The choice of $\BP$ and the construction of the extension of the action have been the source of
developments involving complex reflection groups, their braid groups and Hecke algebras, regular
elements and centralizers, Garside categories and Deligne-Lusztig varieties \cite{Brou3}.

\medskip
With J.~Chuang, we conjecture that the derived equivalence will be perverse, with a non-decreasing
perversity function (for an order as in \S\ref{se:decmat}) \cite{ChRou}. This would imply
the triangularity of the decomposition matrix (cf \S\ref{se:pervexamples}),
a known result (Theorem \ref{th:triangular}).
A conjectural perversity function has been proposed by Craven \cite{Cr}, cf \S\ref{se:disjunction} below.
This implies that the module category of a block with abelian defect groups is determined by
Weyl-group type data, toegther with the perversity function.

\smallskip
When $\ell|(q-1)$, the parabolic subgroup $\BP$ can be chosen to be $F$-stable and the conjecture was
proven by Puig \cite{Pu} (cf \cite[Theorem 23.12]{CaEn2} for a detailed exposition of the principal
block case).  The difficulty is to construct an extension of the action.
As a consequence, for unipotent blocks, the unipotent square part of the decomposition matrix
is diagonal, a fact obtained independently by Hi\ss\ \cite[Korollar 3.2]{Hi}.

\subsubsection{Case of a torus}
\label{se:Brouetori}

We assume that $\ell{\not|}(q-1)$, that $b$ is the principal block idempotent and that $\BL=\BT$ is a
torus. So $D$ is a Sylow $\ell$-subgroup of $G$. Since $C_{\BG}(D)=\BT$, it follows that this torus 
corresponds to a regular $F$-conjugacy class $(w)$ of elements of $W$ (cf \S\ref{se:regular}).
Furthermore, the group
$N_G(D)/C_G(D)$ is isomorphic to $C_W(w\phi)$, a complex reflection group. We denote by $B_d$ its
braid group (cf \S\ref{se:braid}).

Let $w\in W$ such that $(\lambda(w)\phi)^d=\pi \phi^d$.
As explained in \S\ref{se:endomorphisms}, there is a right action of $T\rtimes C_{B_W^+}(\lambda(w)\phi)$ on $Y(\lambda(w))$, hence
a right action on $\Lambda_c(Y(\lambda(w)),\CO)$ commuting with the action of $G$. 

It is conjectured that there is a representative $C$ of $\Lambda_c(Y(\lambda(w)),\CO)$ in the quotient
of the homotopy category of complexes of $\CO(G\times (T\rtimes C_{B_W^+}(\lambda(w)\phi))^\opp)$-modules
by complexes whose restriction to $\CO(G\times T)$ is homotopy equivalent to $0$ with the following
properties:
\begin{itemize}
	\item the right action of $\CO(T\rtimes C_{B_W^+}(\lambda(w)\phi))$ on $C$ factors through an
		action of $\CO N_G(D)b_D$
	\item the resulting complex of $(\CO Gb,\CO N_G(D)b_D)$-bimodules induces a Rickard equivalence
		between the principal blocks of $G$ and $N_G(D)$.
\end{itemize}

This conjecture is known to hold when $X(w)$ is a curve \cite[Corollaire 4.7]{Rou2} and
when $w$ is a Coxeter element \cite{BoRou2,Du1,Du3,DuRou}. In those cases the monoid
$C_{B_W^+}(\lambda(w)\phi)$ is cyclic and its action on $Y(\lambda(w))$ is given by powers of $F$.

\subsubsection{Disjunction of cohomology for Deligne-Lusztig varieties}
\label{se:disjunction}
After extending scalars to $\bar{\BQ}_\ell$, one obtains a version of conjecture \ref{co:Broue2}
that is also an open problem. Restricting to unipotent representations, the crucial missing fact is the
disjunction of the cohomology groups. We state here a conjecture of
\cite{BrouMi2} for the case where $\BL=\BT$ is a torus.

\begin{conj}
	\label{co:disjunction}
	Let $w\in W$ such that $(\lambda(w)\phi)^m=\pi \phi^m$ for some $m\ge 1$. Given $i\neq j$, we have
	$\Hom_{\bar{\BQ}_\ell G}(H^i_c(X(w),\bar{\BQ}_\ell),H^j_c(X(w),\bar{\BQ}_\ell))=0$.
\end{conj}
Craven \cite{Cr} has defined a function $C_m:\Irr_{\bar{\BQ}_\ell}(G,1)\to\BZ$ depending only on the generic
degree of the representation and on $m$ and has conjectured that when $\rho$ occurs in
$H^i_c(X(w),\bar{\BQ}_\ell)$, it occurs in degree $i=C_m(\rho)$.

\smallskip
Conjecture \ref{co:disjunction} (together with Craven's conjecture) is known to hold when $w$ is a Coxeter element \cite{Lu1} and for groups
of rank $2$ \cite{DigMiRou}.
For $\GL_n$, it is known in general (cf \cite{DigMi1} for $m=n-1$ and
\cite[Corollary 3.2]{Du2} and \cite[Theorem 4.3]{BoDuRou} in general).

This conjecture is implied by the refined version of Conjecture \ref{co:Broue2}
discussed in \S \ref{se:Brouetori}.
In the case where $w=w_0$, 
Conjecture \ref{co:disjunction} is older and due to Lusztig
\cite[p.25, line 13]{Lu2}. Also, Conjecture \ref{co:disjunction} was proven earlier
by Lusztig when $w$ is a Coxeter element of minimal length in its class \cite{Lu1}.

\smallskip
Conjecture \ref{co:disjunction} can be extended to $X(b)$, where $b\in B_W^+$ is
such that $(bF)^m=\pi^r F^m$ for some $m,r\ge 1$. The particular case where $b=\pi$ is known to hold
\cite{BoDuRou}.

\smallskip
Conjecture \ref{co:disjunction} can be extended to the case of Deligne-Lusztig varieties associated
to Levi subgroups \cite{DigMi2}.

\smallskip
More recently, Lusztig \cite[\S 7]{Lu4} has conjectured that the disjunction property of
Conjecture \ref{co:disjunction} should hold (in the split case)
for elements $w\in W$ of minimal length in their conjugacy class
and such that the trace of the endomorphism of the Hecke algebra given by $h\mapsto T_w h T_{w^{-1}}$ is in
$\BZ_{\ge 0}[x]$.

\section{Degeneration and genericity}
\label{se:deg}

We consider here the setting of \S\ref{se:Steinberg}.

\subsection{Classifying spaces and character sheaves}
\subsubsection{Completed classifying spaces}

Consider the ring of Witt vectors $R=W(\bar{\BF}_p)$.
Let $\BG_R$ be a reductive algebraic group over $R$ with a maximal torus $\BT_R$ and with an isomorphism 
$\BG_R\times_R \bar{\BF}_p\iso \BG$ restricting to $\BT_R\times_R \bar{\BF}_p\iso \BT_0$. 
We fix an embedding of $R$ into $\BC$ and we denote by
$\BG(\BC)=\BG_R(\BC)$ the associated complex Lie group. We also put $\BT(\BC)=\BT_R(\BC)$.

Specialization provides an isomorphism $\Aut(\BT_R)\iso\Aut(\BT_0)$ and we denote by $\varphi$ the
automorphism of $\BT_R$ lifting $F$. 

We denote by $\varphi$ the automorphism
of $\BT(\BC)$ induced by $F$. The corresponding automorphism
$B\varphi$ of the $\ell$-completed classifying space 
$(B\BT(\BC))_\ell^\wedge$ extends uniquely to an automorphism $\psi$ of
$(B\BG(\BC))_\ell^\wedge$ (\cite[Theorem 1.6]{Fr1} and \cite[Theorem 2.5]{JaMcOl}).

Furthermore, a theorem of Friedlander \cite[Theorem 12.2]{Fr2} (cf also \cite[Theorem 3.1]{BrotMoOl})
shows there is an isomorphism
$$(BG)_\ell^\wedge\iso \bigl((B\BG(\BC))_\ell^\wedge\bigr)^{h\psi}$$
where $h\psi$ denotes taking homotopy fixed points by the group $\BZ$ acting as powers of $\psi$.

\subsubsection{Dependence on $q$}
Consider the group $\Out((B\BG(\BC))_\ell^\wedge)$
of homotopy classes of homotopy automorphisms of $(B\BG(\BC))_\ell^\wedge$. There is an isomorphism
\cite[Theorem 1.2]{AnGro}
$$\Out((B\BG(\BC))_\ell^\wedge)\iso N_{\GL(Y(\BT_0)\otimes\BZ_\ell)}(W,
\{\BZ_\ell \beta\}_{\beta\in\Phi^\vee})/W.$$

\smallskip
Given $\alpha\in\Out((B\BG(\BC))_\ell^\wedge)$, the space
$(B\BG(\BC))_\ell^\wedge)^{h\alpha}$ depends only on the closed subgroup
$\overline{\langle\alpha\rangle}$ of $\Out((B\BG(\BC))_\ell^\wedge)$ \cite[Corollary 2.5]{BrotMoOl}, where
we use the $\ell$-adic topology.

\medskip
	When $(\BG,F)$ is split, then $\varphi$ is the automorphism $x\mapsto x^q$ and 
	$\psi$ is the unstable Adams operation $\psi^q$. Note that the unstable Adams operation
	$\psi^q$ is defined more generally for $q\in \BZ_\ell^\times$.

	In general, when $\BG$ is simple and $F$ is a Frobenius endomorphism,
	then the element of $N_{\GL(Y(\BT_0)\otimes\BZ_\ell)}(W,
	\{\BZ_\ell \beta\}_{\beta\in\Phi^\vee})/W$ induced by $\psi$ is of the form $\sigma\cdot (q\id)$
	where $\sigma$ has finite order
	and $\psi=\sigma\psi^q$ (up to homotopy). In types ${^2A}_n$, ${^2D}_{2n+1}$ and
	${^2E_6}$, one has also $\psi=\psi^{-q}$.

	The description $\psi=\sigma\cdot (q\id)$
	still works for types ${^2B}_2$ and ${^2F}_4$ (resp. ${^2G}_2$) when $2$ 
	(resp. $3$) is a square modulo $\ell$.

	\medskip
	Assume $\ell$ is odd.
	The space $(BG)_\ell^\wedge$ depends only on the order $d$ of $q$ in $\BF_\ell^\times$ and on 
	$\nu_\ell(q^d-1)$ \cite[\S 3 and Proposition 3.2]{BrotMoOl}.
	
	Also, $B({^2A_n}(q))_\ell^\wedge\simeq B(A_n(q'))_\ell^\wedge$,
	$B({^2D_{2n+1}}(q))_\ell^\wedge\simeq B(D_{2n+1}(q'))_\ell^\wedge$ and
	 $B({^2E_6}(q))_\ell^\wedge\simeq B(E_6(q'))_\ell^\wedge$ if
	$q$ and $-q'$ have the same order $d$ in $\BF_\ell^\times$ and 
	$\nu_\ell(q^d-1)=\nu_\ell((-q')^d-1)$ \cite[Proposition 3.3]{BrotMoOl}.

	Note that $(BG)_\ell^\wedge$ determines the thick subcategory of
	$D^b(\BF_\ell G)$ generated by the trivial module, so this triangulated category
	has the same genericity properties.

\subsubsection{Classifying spaces of loop groups}
\label{se:Bloop}

We assume $(\BG,F)$ is split and put $\eps=1$ or 
$(\BG,F)$ has type ${^2A}_n$, ${^2D}_{2n+1}$ or ${^2E_6}$ and put $\eps=-1$.
We assume $\ell|\eps q-1$.

\smallskip
We have a family of spaces
$(B\BG(\BC)^\wedge_\ell)^{h\psi_{1+\ell\hbar}}$ over $\BZ_\ell$,
constant over $\BZ_\ell^\times$-orbits:
$$\xymatrix{
(B\BG(\BC)^\wedge_\ell)^{h\psi_{1+\ell\hbar}} \ar[d] \\
\BZ_{\ell}=\{\hbar\}\ar[r] \ar[dr]_{\nu_\ell}& \BZ_\ell/\BZ_\ell^\times\ar[d]_\sim \\
& \BZ_{\ge 0}\cup\{\infty\}
}$$

One could attempt to make sense of this as a continuous family, and then make sense of
the limit of those spaces as $\hbar\to 0$ to obtain
$$"\lim_{\hbar\to 0}"(B\BG(\BC)^\wedge_\ell)^{h\psi_{1+\ell\hbar}}\simeq
(B\BG(\BC)_\ell^\wedge)^{h\mathrm{id}}=L(B\BG(\BC))^\wedge_\ell\simeq
B(L\BG(\BC))^\wedge_\ell,$$
where $L(X)=\mathrm{Maps}(S^1,X)$ is the free loop space and in
particular $L\BG(\BC)$ is the loop group associated to $\BG(\BC)$.
So, from the point of view of $\ell$-completed classifying spaces, the loop group $L\BG(\BC)$ appears as
$\BG(\BF_1)$.

In other terms,
$$"\lim_{\nu_l(\eps q-1)\to \infty}"B\BG(\BF_q)_\ell^\wedge\simeq
B(L\BG(\BC))^\wedge_\ell.$$

So the space $B(L\BG(\BC))^\wedge_\ell$ appears as a degeneration of the family of spaces
$B\BG(\BF_q)^\wedge_\ell$ for varying $q$. Here, we use the abusive notation $\BG(\BF_q)$ to denote
a possibly twisted group in a family.

As a consequence, we have also a description of the limit ("generic version") of the thick subcategory
of $D^b(\BF_\ell G)$ generated by the trivial module as $\nu_\ell(\eps q-1)\to \infty$: it is the homotopy
category of perfect $A_\infty$-modules over the $A_\infty$-algebra $H^*(BL\BG(\BC),\BF_\ell)$, since
the thick subcategory of $D^b(\BF_\ell G)$ generated by the trivial module is equivalent to perfect
$A_\infty$-modules over $\Ext^*_{\BF_\ell G}(\BF_{\ell},\BF_{\ell})$.

Note that while the family of algebras $H^*(G,\BF_\ell)$ stabilizes \cite[Theorem 18]{KiKo},
the stabilization does not hold when the $A_\infty$-algebra structure is taken into account (cf
Remark \ref{re:Gm} below).

\begin{rem}
	\label{re:Gm}
When $\BG=\BG_m$, we have $BL\BG(\BC)=BL\BC^\times\simeq S^1\times\BC\BP^\infty$, a space whose
mod-$\ell$ cohomology is formal as an algebra.
The
$A_\infty$-structure on $H^*(B\BF_q^\times,\BF_\ell)$ can be chosen so that there is a single
	higher multiplication $m_{\ell^r}$, where $r=\nu_\ell(q-1)$ \cite[Appendix B, Example 2.2]{Mad}.
	This higher multiplication disappears in
the limit $r\to\infty$.
\end{rem}

\begin{rem}
The relation between the cohomology of the finite group and that of the loop
	group becomes much more subtle for small $\ell$ (and $r$). An approach using the string topology
	is given in \cite{GroLa}.
\end{rem}

\begin{rem}
	Considering the usual topology instead of the $\ell$-adic one, we obtain
	$\lim B\BG(\BF_q)_\ell^\wedge=B\BG(\bar{\BF}_q)_\ell^\wedge\simeq B\BG(\BC)_\ell^\wedge$
	\cite[Theorem 1.4]{FrMi}.
\end{rem}

\subsubsection{Rigidification and character sheaves}
\label{se:rigidification}
In the previous section, we explained how to obtain, under some assumptions,
a generic version of the modular representation theory of $G$, in the neighborhood of the
trivial representation.

To move away from the neighborhood of the trivial representation, we consider a more rigid version of
$B(L\BG(\BC))_\ell^\wedge$.
There is a homotopy equivalence $B(L\BG(\BC))\simeq \frac{\BG(\BC)}{h\BG(\BC)}$, the homotopy adjoint
quotient.

We can now consider the derived category of 
$\CD$-modules on the stack $\frac{\BG(\BC)}{\BG(\BC)}$, or constructible sheaves with
$k$-coefficients. This is the $\BG(\BC)$-equivariant derived category of $\BG(\BC)$,
for the adjoint action, and it has a thick subcategory of unipotent objects,
also called the derived unipotent character sheaves, providing a non semi-simple enrichment of Lusztig's
theory \cite[Definition 6.8]{BeZNa}. It is conjectured that
the principal series part of this triangulated category (i.e., its principal block)
is equivalent, for $\ell$ not too small, to the derived category of differential
graded modules over $k[\Gh\times\Gh^*]\rtimes W$, where
$\Gh$ is the Lie algebra of $\BT_0$ over $k$.
A similar result is known for the adjoint quotient
$\frac{\mathrm{Lie}\BG}{\BG}$ \cite{Rid}.

For $\BG=\GL_n$, all unipotent character sheaves are in the principal series,
so the description coincides with our degeneration approach in \S\ref{se:dege} below.

\smallskip
An important problem is to find a conjecture for an algebraic description of the derived unipotent
character sheaves, beyond the principal series, starting with the case $\BG=\mathrm{Sp}_4$.
This is related to the problem of finding a canonical generic description of the category of
unipotent representations in characteristic zero, cf \cite{Lu3}.

\smallskip

Note that the category breaks down according to Harish-Chandra series,
but it would be desirable to find a description that does not use
cuspidal objects.

\subsubsection{General $d$}
The first constructions of \S\ref{se:Bloop} can be performed without the assumption that $\ell|\eps q-1$.
Denote by $d$ the order of $\eps q$ in $\BF_\ell^\times$. Let $\zeta$ be a primitive $d$-th root
of unity in $\BZ_\ell$.

One can consider the family of spaces
$(B\BG(\BC)^\wedge_\ell)^{h\psi_{\zeta+\ell\hbar}}$ over $\BZ_\ell=\{\hbar\}$ and
$$"\lim_{\nu_l(\Phi_d(\eps q))\to \infty}"B\BG(\BF_q)_\ell^\wedge\simeq
"\lim_{\hbar\to 0}"(B\BG(\BC)^\wedge_\ell)^{h\psi_{\zeta+\ell\hbar}}\simeq
L(B\BG(\BC)_\ell^\wedge)^{h\mu_d}.$$
Here $\mu_d$ is the cyclic group of order $d$ acting on $B\BG(\BC)_\ell^\wedge$ by $\psi_x$, 
$x\in\mu_d(\BZ_\ell)$.

\smallskip
The space $(BG(\BC)^\wedge_\ell)^{h\mu_d}$ is an $\ell$-compact group \cite{Gro}.
Its "Weyl group" is a complex (or rather $\ell$-adic) reflection group,
not a Coxeter group in general. To proceed as in \S\ref{se:rigidification} we would need an
appropriately rigidified version of the space $L(BG(\BC)^\wedge_\ell)^{h\mu_d}$.

\begin{rem}
	In \cite{KeMalSe}, Kessar, Malle and Semeraro explain how to understand Alperin's conjecture
in the setting of $\ell$-completed classifying spaces. One can expect there is a framework which
encompasses both the cohomological aspects, which was our starting point, and the character counts, which
they study.
\end{rem}

\subsection{Degeneration}
\label{se:dege}

\subsubsection{Degeneration of group algebras of abelian $\ell$-groups}
\label{se:degenerateabelian}
Let $P$ be an abelian $\ell$-group isomorphic to $(\BZ/\ell^r)^n$. Let
$V=J(\BF_{\ell}P)/J(\BF_{\ell}P)^2$. This is an $n$-dimensional vector space over $\BF_{\ell}$.
Fix a morphism of $\BF_{\ell}$-modules $\sigma:V\to J(\BF_{\ell}P)$ that is a right inverse to the quotient map
$J(\BF_{\ell}P)\to V$. The map $\sigma$ extends uniquely to a morphism of $\BF_{\ell}$-algebras
$S(V)\to \BF_{\ell}P$. That morphism induces an isomorphism
\begin{equation}
	\label{eq:groupalgasS}
	S(V)/(v^{\ell^r})_{v\in V}\iso \BF_{\ell}P.
\end{equation}

\smallskip
Consider now a finite $\ell'$-group $E$ acting on $P$. The vector space $V$ is an $\BF_{\ell}E$-module.
Since $J(\BF_{\ell}P)$ is a semi-simple
$\BF_{\ell}E$-module, there exists a $\sigma$ as above that is a morphism of $\BF_{\ell}E$-modules. The
	isomorphism (\ref{eq:groupalgasS}) is equivariant for the action of $E$, hence it
	extends to an isomorphism of $\BF_\ell$-algebras
$$\bigl(S(V)/(v^{\ell^r})_{v\in V}\bigr)\rtimes E\iso \BF_{\ell}(P\rtimes E).$$

\medskip
Consider now a general finite abelian $\ell$-group $P$ acted on by a finite $\ell'$-group $E$.
There exists an $E$-stable decomposition $P=P_1\times\cdots\times P_m$ such that 
$P_i\simeq(\BZ/\ell^{r_i})^{n_i}$ for some $r_i$ and $n_i$. Put $V_i=J(\BF_{\ell}P_i)/J(\BF_{\ell}P_i)^2$.
The construction above provides
an isomorphism of $\BF_\ell$-algebras
$$\bigl(S(V)/(\bigcup_i \{v^{\ell^{r_i}}\}_{v\in V_i})\bigr)\rtimes E\iso \BF_{\ell}(P\rtimes E).$$

\medskip
Consider the graded $\BF_{\ell}[t]$-algebra $A=\BF_{\ell}[t]\otimes\Lambda(V)\otimes S(V)$, where 
$\BF_{\ell}[t]\otimes \BF_{\ell}\otimes S(V)$ is in degree $0$ and $\BF_{\ell}\otimes V\otimes \BF_{\ell}$
is in degree $-1$.

We define a structure of differential $(\BF_{\ell}[t]\otimes\BF_{\ell}\otimes S(V))$-algebra on
$A$ by setting
$d(1\otimes v\otimes 1)=t\otimes 1\otimes v^{\ell^{r_i}}$ for $v\in V_i$.

We have $H^i(\BF_{\ell}(t)\otimes_{\BF_{\ell}[t]}A)=0$ for $m\neq 0$ and 
$$H^0(\BF_{\ell}(t)\otimes_{\BF_{\ell}[t]}A)=(\BF_{\ell}(t)\otimes S(V))/(\bigcup_i 
\{tv^{\ell^{r_i}}\}_{v\in V_i}\simeq
\BF_{\ell}(t)\otimes S(V)/(\bigcup_i \{v^{\ell^{r_i}}\}_{v\in V_i}).$$

So, the algebra $\BF_{\ell}(P\rtimes E)$ is, up to quasi-isomorphism, a deformation of the graded
algebra $\bigl(\Lambda(V)\otimes S(V)\bigr)\rtimes E$. The derived category of
$\BF_{\ell}(P\rtimes E)$-modules is a deformation of the derived category of dg modules over the graded
algebra (with zero differential) $\bigl(\Lambda(V)\otimes S(V)\bigr)\rtimes E$.

\medskip
Koszul duality provides an equivalence from the derived category of finitely generated
differential graded modules over the graded
algebra (with zero differential) $\bigl(\Lambda(V)\otimes S(V)\bigr)\rtimes E$ to
the derived category of finitely generated
differential graded modules over the graded algebra $S(V^*\oplus V)\rtimes E$
(here $V$ is in degree $0$ and $V^*$ in degree $2$).

\medskip
To summarize, $D^b(\BF_{\ell}(P\rtimes E))$ degenerates
into the derived category of differential graded coherent sheaves on the orbifold $[(V\times V^*)/E]$.

\subsubsection{Genericity of perverse equivalences}
\label{se:genericperverse}

The discussion here is based on joint work with David Craven \cite{CrRou2}.
We consider the setting of \S\ref{se:BroueconjP} and we assume to simplify that $b$ is the principal
block. So $D$ is a Sylow $\ell$-subgroup of $G$ and 
there is an isomorphism of algebras $kN_G(D)b_D\simeq kD\rtimes E$ where
$E=N_G(D)/C_G(D)$.

It is conjectured that there is a perverse equivalence between $kD\rtimes E$ and 
$kGb$, with a specific perversity function $\pi:\Irr_{\BC}(E)\iso\Irr_k(E)\to\BZ$. That function depends
only on the type of the group $G$ and on $d$, not on $q$ or $\ell$.

As explained in \S\ref{se:localnondef}, the group $E$ is a reflection group. We denote by $K_E$ the field
of definition
of its reflection representation $V$ and by $\CO_E$ the ring of integers of $K_E$. Let 
$R=\CO_E[|W|^{-1}]$ and let $V_R$ be an $RE$-module, finitely generated and projective over $R$, such that 
$V\simeq K\otimes_R V_R$.

We conjecture that the function $\pi$ defines a $t$-structure on the derived category of differential
graded modules
over the graded algebra $(\Lambda(V_R)\otimes S(V_R))\rtimes E$, where $(R\otimes S(V_R))\rtimes E$ is in degree
$0$ and $V_R\otimes R$ in degree $-1$.
The heart $\CA$ of that $t$-structure would be 
a ``generic version'' of $kGb$, i.e., a limit as $\nu_\ell(\Phi_d(\eps q))\to \infty$. 
A ridigity property of perverse simple objects would show that the classes of
the indecomposable projective objects of $\CA\otimes_R K_E$ expressed in terms of the classes
of the simple $K_E E$-modules would give the transpose of the square unipotent part of the decomposition 
matrix of the principal $\ell$-block of $G$ for $\ell$ large enough. Note that the presence of a double
grading on $(\Lambda(V)\otimes S(V))\rtimes E$ leads to a two-variable deformation of the matrix.

\begin{rem}
	The discussion generalizes to the case of non-principal blocks. The block
	of the normalizer is isomorphic to a twist of the group algebra of the semi-direct product by
	a $2$-cocycle but that cocycle is expected to be always trivial.
\end{rem}

\begin{rem}
We expect the algebra $S(V\oplus V^*)\rtimes E$ to control generic
aspects of the modular representation theory of $G$. This algebra admits deformations as rational
Cherednik algebras and the "$t=0$" case is expected to relate to unipotent representations
	\cite{BoRou4,Bo}.
\end{rem}

\subsubsection{Hilbert schemes}
We discuss here joint work with Olivier Dudas \cite{DuRou}.

Assume that $\BG=\GL_n$. Let $m=\lfloor \frac{n}{d}\rfloor$. We have $V\simeq K_E^m$ and
$W\simeq (\BZ/d)^m\rtimes\GS_m$. Let $X_d$ be the minimal resolution of $\BA_{K_E}^2/(\BZ/d)$, where
$\BZ/d$ is embedded in $\SL_2(K_E)$.

\smallskip
Let $\Hilb^m(X_d)$ be the Hilbert scheme of $m$ points on $X_d$ and
$\pi:\Hilb^n(X_d)\to S^m(X_d)$ be the Hilbert-Chow map. Let $f:\BA^{2m}\to S^m(X_d)$ be the
quotient map by $(\BZ/d)^m\rtimes\GS_m$.

Combining Koszul duality with the derived McKay equivalence,
we obtain an equivalence between the derived category of differential
graded $(\Lambda(V)\otimes S(V))\rtimes E$-modules and the derived category of dg coherent sheaves 
on $\Hilb^m(X_d)$, where we consider the $\BG_m$-action on
$X_d$ coming from its action on $\BA^2$ with weights $0$ and $-2$. The conjecture in 
\S\ref{se:genericperverse}
implies the existence of a particular $t$-structure on that derived category.

\smallskip

When $d=1$ and $\eps=-1$, the combinatorics of Macdonald polynomials can be used to obtain a conjectural
combinatorial formula for the two-parameter deformed decomposition matrix of $U_n(q)$.
That conjecture has been checked for $n\le 11$, using the determination of the decomposition matrices
in \cite{DuMal1}.

\section{Appendix}
\subsection{Representations}

\subsubsection{Categories}
Let $A$ be an algebra over a commutative regular local
noetherian ring $R$ and assume $A$ is a free $R$-module of finite rank.
By module, we mean left module. We identify right $A$-modules with left modules
for the opposite algebra $A^\opp$.

Given $M$ an $A$-module, we put $M^*=\Hom_R(M,R)$, a right $A$-module.

We denote by $\Irr(A)$ the set of isomorphism classes of simple $A$-modules.

\smallskip
We denote by $A\mMod$ the abelian category of finitely generated $A$-modules.
We denote by $G_0(A)$ the Grothendieck group of $A\mMod$.

We denote by $D^b(A)$ (resp. $\Ho^b(A)$) the derived (resp. homotopy)
category of bounded complexes of finitely generated $A$-modules. We denote by $A\mperf$ the
full subcategory of $D^b(A)$ of complexes quasi-isomorphic to bounded complexes of finitely
generated projective $A$-modules.

Let $A\mstab$ be the triangulated category quotient $D^b(A)/A\mperf$.
When $A$ is symmetric as an $R$-algebra, the inclusion $A\mMod\hookrightarrow D^b(A)$ induces
an equivalence of categories from 
the additive category quotient of $A\mMod$ by its subcategory of $A$-modules
of the form $A\otimes_R V$ where $V\in R\mMod$, to $A\mstab$.

\subsubsection{Equivalences}
\label{se:equivalences}

Let $A$ and $B$ be two finite-dimensional algebras over a commutative noetherian ring $R$.

Let $C$ be a bounded complex of $(A,B)$-bimodules, all of whose terms are finitely generated
and projective as left $A$-modules and as right $B$-modules. Assume there is a complex $L$ (resp. $M$)
of $(A,A)$-bimodules (resp. $(B,B)$-bimodules) such that there are isomorphisms of 
complexes of $(A,A)$-bimodules and $(B,B)$-bimodules
$$C\otimes_BC^*\simeq A\oplus L \text{ and }C^*\otimes_A C\simeq B\oplus M.$$
We say that $M$ induces a
\begin{itemize}
	\item Morita equivalence if $C^i=0$ for $i\neq 0$ and $L=M=0$
	\item Rickard equivalence if $L$ and $M$ are homotopy equivalent to $0$
	\item derived equivalence if $L$ and $M$ are acyclic
	\item stable equivalence if $L$ and $M$ are perfect.
\end{itemize}

These conditions ensure that $C\otimes_B-$ induces an equivalence

\begin{itemize}
	\item (Morita) $B\mMod\iso A\mMod$
	\item (Rickard) $\Ho^b(B)\iso \Ho^b(A)$
	\item (derived) $D^b(B)\iso D^b(A)$
	\item (stable) $B\mstab\iso A\mstab$
\end{itemize}

\subsubsection{Finite groups}
\label{se:appfinite}
Let $G$ be a finite group. We put $\Irr_R(G)=\Irr(RG)$.
Consider a prime $p$ and a finite field extension $K$ of $\BQ_p$. Let $\CO$ be its ring of integers and $k$
the residue field. We assume that $K$ contains all $|G|$-th roots of unity. This ensures that
$KG$ is a product of matrix algebras over $K$ and that all simple $kG$ modules are absolutely simple.

\smallskip
Let $M$ be a finitely generated $KG$-module. There exists an $\CO G$-module $M'$ that is free
over $\CO$ and such that $M'\otimes_{\CO}K\simeq M$. Let $M''=M'\otimes_{\CO}k$. The class $[M'']$
in $G_0(kG)$ depends only on $[M]\in G_0(KG)$ and we put $\mathrm{dec}([M])=[M'']$. This defines a morphism
of abelian groups, the {\em decomposition map}, $\mathrm{dec}:G_0(KG)\to G_0(kG)$. The {\em decomposition matrix}
is the matrix of $\mathrm{dec}$ in the bases $\Irr_k(G)$ (columns) and $\Irr_K(G)$ (rows).

\subsection{Braid groups and Hecke algebras}
\subsubsection{Braid groups}
\label{se:braid}
Let $V$ be a finite dimensional complex vector space.
A {\em reflection} $s$ of $V$ is a finite order automorphism of $V$ such that $\ker(s-1)$ is
a hyperplane.

Let $W$ be a finite subgroup of $\GL(V)$ generated by reflections (a {\em complex reflection group}).
Let $\CR$ be the set of reflections in $W$
and $\CA=\{\ker(s-1)\}_{s\in\CR}$ be the set of reflecting hyperplanes.

We put $V^{reg}=V\setminus\bigcup_{H\in\CA}H$. The group $W$ acts
freely on $V^{reg}$,
i.e., the quotient map $q:V^{reg}\to V^{reg}/W$ is unramified.

Let 
$x_0 \in V^{reg}$. The {\em braid group} of $W$ is $B_W=\pi_1(V^{reg}/W,q(x_0))$.
The map $q$ gives a bijection from (homotopy classes of) paths in $V^{reg}$ starting at $x_0$ and ending in
$W(x_0)$ to (homotopy classes of) loops in $V^{reg}/W$ based at $q(x_0)$, and we will identify those
two types of objects.
There is a surjective  morphism $B_W\to W$: it sends $w$ to the homotopy class of
a path in $V^{reg}$ from $x_0$ to $w(x_0)$.
We denote by $\Bpi\in B_W$ the homotopy class of the path $t\mapsto \exp(2i\pi t)x_0$. This is a
central element of $B_W$.

\subsubsection{Hecke algebras}
\label{se:Hecke}
Given $H\in\CA$, let $e_H$ be the order of the fixator of $H$ in $W$.
Let $R=\BZ[\{q_{H,r}^{\pm 1}\}_{H\in\CA/W, 0\le r<e_H}]$.

We define the {\em Hecke algebra} $\CH=\CH(W)$ of $W$
as the quotient of the group algebra $RB_W$ by the ideal generated by
$\prod_{0\le r<e_H}(\sigma_H-q_{H,r})$, where $H$ runs over $\CA$ and $\sigma_H$ is a generator
of the monodromy around the image of $H$ in $V/W$
\cite[Definition 4.21]{BrouMalRou}.

\smallskip
The specialization $q_{H,r}\mapsto\exp(2i\pi r/e_H)$ of $\CH$ is the group algebra $\BZ W$.

\subsubsection{Regular elements}
\label{se:regular}
We recall some constructions and results of Springer \cite{Sp}.

Let $\sigma$ be an element of finite order of $N_{\GL(V)}(W)$. Let
$w\in W$ and let $v\in V^{reg}$ be an eigenvector of $w\sigma$ with eigenvalue $\zeta$.  Let
$d$ be the order of $\zeta$.
The element $w\sigma$ is said to be {\em $\zeta$-regular}, or {\em $d$-regular}. If $w'\in W$ and $w'\sigma$ is
$\zeta$-regular, then $w'\sigma$ is $W$-conjugate to $w\sigma$.

\smallskip
Let $V_\zeta=\ker(w\sigma-\zeta)$. The group $C_W(w\sigma)$ acting on $V_\zeta$ is a reflection group.

The inclusion $V_\zeta\hookrightarrow V$ induces an
isomorphism $\iota_\zeta:V_{\zeta}/C_W(w\sigma)\iso (V/W)^{\mu_d}$, where
$\mu_d=\{\zeta^n\id_V\}_{n\in\BZ/d}$.

Assume $\zeta=\exp(2i\pi/d)$ and $x_0=v$. There exists 
$\Bw_d\in B_W$ such that $(\Bw_d\sigma)^d=\Bpi\sigma^d$ \cite[Proposition 6.5]{BrouMi2}.
When $\sigma=1$ we can take for $\Bw_d\in B_W$ the homotopy class of the path
$t\mapsto\exp(2i\pi t/d)x_0$.

The map $\iota_\zeta$ induces a morphism $B_{C_W(w\sigma)}=\pi_1(V_{\zeta}^{reg}/C_W(w\sigma),q(x_0))\to
B_W=\pi_1(V^{reg}/W,q(x_0))$. Its image is contained in $C_{B_W}(\Bw_d\sigma)$.

\subsubsection{Real reflection groups}
\label{se:real}

We assume now that $V=V_\BR\otimes_{\BR}\BC$ and $W$ is a subgroup of $\GL(V_{\BR})$.
All reflections of $W$ have order $2$.

Fix a connected component $C$ of the space $V_{\BR}\cap V^{reg}$ and let $\bar{C}$ be its closure.
Let $S$ be the subset of
$\CR$ of reflections $s$ such that $\ker(s-1)\cap \bar{C}$ has codimension $1$ in $V_{\BR}$. Then
$(W,S)$ is a Coxeter group. We denote by $l:W\to\BZ_{\ge 0}$ its length function: given $w\in W$,
the integer $l(w)$ is the minimal $m$ such that $w=s_{i_1}\cdots s_{i_m}$ for some
$s_{i_1},\ldots,s_{i_m}\in S$.

\smallskip
Choose now $x_0\in C$. Given $s\in S$, let $\sigma_s\in B_W$ be the homotopy class of the path that is
the concatenation of $t\mapsto x_0+tix_0$, $t\mapsto (1-t)x_0+ts(x_0)+ix_0$ and $t\mapsto s(x_0)+(1-t)ix_0$.

There is an isomorphism
\begin{equation}
	\label{eq:isobraid}
\langle (b_s)_{s\in S}\ |\ \underbrace{b_sb_tb_s\cdots}_{m_{st}\ \text{terms}}=
\underbrace{b_tb_sb_t\cdots}_{m_{st}\ \text{terms}} \forall s,t\in S\}\iso
B_W,\ b_s\mapsto\sigma_s
\end{equation}
where $m_{st}$ is the order of $st$ \cite{Bri}. We identify $B_W$ with the group on the left
side of (\ref{eq:isobraid}) and we denote by $B_W^+$ its submonoid generated by $(b_s)_{s\in S}$.

There is a map $\lambda:W\to B_W$ given by $\lambda(w)=b_{s_1}\cdots b_{s_r}$ if
$w=s_1\cdots s_r$ is any reduced decomposition of $w\in W$ with $s_i\in S$.
Denote by $w_0$ the longest element of $W$. We have $\pi=\lambda(w_0)^2$.

\medskip
Let $x$ be an indeterminate and let $H(W)$ be the "usual" Hecke algebra of $W$, i.e., the
$\BZ[x^{\pm 1}]$-algebra generated by $(T_s)_{s\in S}$ with relations
$$(T_s-x)(T_s+1)=0,\ \underbrace{T_sT_tT_s\cdots}_{m_{st}\ \text{terms}}=
\underbrace{T_tT_sT_t\cdots}_{m_{st}\ \text{terms}} \text{ for } s,t\in S.$$

The isomorphism (\ref{eq:isobraid}) induces an isomorphism between $H(W)$ and
the specialization of $\CH(W)$ at $q_{H,0}\mapsto x,\ q_{H,1}\mapsto -1$.

\end{document}